\documentclass[pdflatex,sn-mathphys-num]{sn-jnl}
\usepackage{amsmath,amssymb,amsfonts,amsthm}%
\usepackage{color}
\allowdisplaybreaks

\theoremstyle{thmstyleone}%
\newtheorem{theorem}{Theorem}[section]
\newtheorem{lemma}[theorem]{Lemma}%
 \newtheorem{corollary}[theorem]{Corollary}
\theoremstyle{thmstyletwo}%
\theoremstyle{thmstylethree}%
\numberwithin{equation}{section}

\begin{document}
\title[Gohberg-Krupnik Localisation for Discrete Wiener-Hopf Operators]
{{Gohberg}-Krupnik Localisation {for} 
Discrete Wiener-Hopf Operators on Orlicz Sequence Spaces}

\author*[1]{\fnm{Oleksiy} \sur{Karlovych}}
\email{oyk@fct.unl.pt}
\author[1]{\fnm{Sandra Mary} \sur{Thampi}}
\email{s.thampi@campus.fct.unl.pt}
\equalcont{{These authors contributed equally to this work.}\\ \\
To Professor Sergei Grudsky on the occasion of his 70th birthday}
\affil*[1]{
\orgdiv{Centro de Matem\'atica e Aplica\c{c}\~oes, Departamento de Matem\'atica}, 
\orgname{Faculdade de Ci\^encias e Tecnologia, Universidade Nova de Lisboa}, 
\orgaddress{
\street{Quinta da Torre}, 
\city{Caparica}, 
\postcode{2829--516}, 
\country{Portugal}}}
\abstract{Let $\Phi$ be an $N$-function whose Matuszewska-Orlicz indices
satisfy $1<\alpha_\Phi\le\beta_\Phi<\infty$. Using these indices, we 
introduce ``interpolation friendly" classes of Fourier multipliers 
$M_{[\Phi]}$ and $M_{\langle\Phi\rangle}$ such that 
$M_{[\Phi]}\subset M_{\langle\Phi\rangle}\subset M_\Phi$, where
$M_\Phi$ is the Banach algebra
of all Fourier multipliers on the reflexive Orlicz sequence space 
$\ell^\Phi(\mathbb{Z})$. 
Applying the Gohberg-Krupnik localisation in the corresponding Calkin
algebra, the study of Fredholmness of the discrete Wiener-Hopf
operator $T(a)$ with $a\in M_{\langle\Phi\rangle}$ is reduced
to that of $T(a_\tau)$ for certain, potentially easier to study, local
representatives $a_\tau\in M_{[\Phi]}$ of $a$ at all points
$\tau\in[-\pi,\pi)$. 
}
\keywords{%
Discrete Wiener-Hopf operator,
Fredholmness,  
Gohberg-Krupnik local principle,
Orlicz sequence space,
Matuszewska-Orlicz indices.}
\maketitle
\section{Introduction and the main result}
For a Banach space $\mathcal{X}$, let $\mathcal{B}(\mathcal{X})$ and
$\mathcal{K}(\mathcal{X})$ denote the Banach algebra of all bounded linear
operators on $\mathcal{X}$ and its closed two-sided ideal consisting of 
all compact linear operators on $\mathcal{X}$, respectively. An operator
$A\in\mathcal{B}(\mathcal{X})$ is said to be Fredholm on $\mathcal{X}$
if $\dim\operatorname{Ker}A<\infty$ and 
$\dim\mathcal{X}/\operatorname{Im}A<\infty$, where 
$\operatorname{Ker}A:=\{x\in\mathcal{X}:Ax=0\}$ and 
$\operatorname{Im}A:=A(\mathcal{X})$ are the kernel and the range 
of $A$, respectively. We will abbreviate the coset 
$A+\mathcal{K}(\mathcal{X})$ of the Calkin algebra 
$\mathcal{B}(\mathcal{X})/\mathcal{K}(\mathcal{X})$ to $A^\pi$.

A function $\Phi:[0,\infty)\to[0,\infty)$ is said to be an $N$-function
if it is continuous, increasing, convex, $\Phi(t)=0$ if and only if $t=0$,
and
\[
\lim_{t\to 0^+}\frac{\Phi(t)}{t}=0,
\quad
\lim_{t\to\infty}\frac{\Phi(t)}{t}=\infty.
\]
For properties of $N$-functions, see, e.g., \cite[Ch.~1, \S~1]{KR61} or \cite[Ch.~8]{M89}.
Let $\ell^0(\mathbb{Z})$ be the linear space of all sequences 
$f:\mathbb{Z}\to\mathbb{C}$. For an $N$-function $\Phi$, the Orlicz sequence
space $\ell^\Phi(\mathbb{Z})$ consists of all sequences 
$f=\{f_n\}_{n\in\mathbb{Z}}\in\ell^0(\mathbb{Z})$ such that
\[
\varrho_\Phi(f/\lambda):=\sum_{n\in\mathbb{Z}}\Phi(|f_n/\lambda|)<\infty
\]
for some $\lambda=\lambda(f)>0$. It is well known that $\ell^\Phi(\mathbb{Z})$
is a rearrangement-invariant Banach sequence space with respect to the 
Luxemburg-Nakano norm given by
\[
\|f\|_{\ell^\Phi(\mathbb{Z})}:=\inf\{\lambda>0:\varrho_\Phi(f/\lambda)\le 1\}
\]
(cf., e.g., \cite[Ch.~4, Theorem~8.9]{BS88}). It is easy to see that 
if $1<p<\infty$, then $\Phi_p(t)=t^p$ is an $N$-function and 
$\ell^{\Phi_p}(\mathbb{Z})$ is nothing but the standard Lebesgue sequence 
space $\ell^p(\mathbb{Z})$.

Let $\mathcal{P}'$ be the space of $2\pi$-periodic distributions (see, e.g.,
\cite[Ch.~3 and~5]{B73})
and let 
$S_0(\mathbb{Z})$ denote the set of all finitely supported sequences.
For $a\in\mathcal{P}'$ and $\varphi\in S_0(\mathbb{Z})$, we define
the convolution $a*\varphi$ as the sequence
\[
(a*\varphi)_j=\sum_{k\in\mathbb{Z}}\widehat{a}_{j-k}\varphi_k,
\quad
j\in\mathbb{Z},
\]
where $\{\widehat{a}_j\}_{j\in\mathbb{Z}}$ is the sequence of the Fourier
coefficients of the distribution $a$. By $M_\Phi$ we denote
the collection of all distributions $a\in\mathcal{P}'$ for which 
$a*\varphi\in \ell^\Phi(\mathbb{Z})$ whenever $\varphi\in S_0(\mathbb{Z})$ 
and
\begin{equation}\label{eq:multiplier-norm}
\|a\|_{M_\Phi}:=\sup\left\{
\frac{\|a*\varphi\|_{\ell^\Phi(\mathbb{Z})}}
{\|\varphi\|_{\ell^\Phi(\mathbb{Z})}}
\ :\
\varphi\in S_0(\mathbb{Z}), \ \varphi\ne 0
\right\}<\infty.
\end{equation}
Suppose the space $\ell^\Phi(\mathbb{Z})$ is separable. Then $S_0(\mathbb{Z})$ 
is dense in $\ell^\Phi(\mathbb{Z})$ (see, e.g., \cite[Lemma~3]{KT24}). 
In this case, for $a\in M_\Phi$, the operator
from $S_0(\mathbb{Z})$ to $\ell^\Phi(\mathbb{Z})$ defined by 
$\varphi\mapsto a*\varphi$ extends to a bounded operator 
\begin{equation}\label{eq:Laurent}
L(a):\ell^\Phi(\mathbb{Z})\to \ell^\Phi(\mathbb{Z}),
\quad
\varphi\mapsto a*\varphi,
\end{equation}
which is referred to as the Laurent operator with symbol $a$. 
If, in addition, $\ell^\Phi(\mathbb{Z})$ is reflexive, then
\begin{equation}\label{eq:embedding-of-multipliers}
\|a\|_{L^\infty(-\pi,\pi)}\le \|a\|_{M_\Phi}
\quad\mbox{for all}\quad
a\in M_\Phi.
\end{equation}
Here and in what follows, $L^\infty(-\pi,\pi)$ denotes the Banach space 
of all essentially bounded $2\pi$-periodic functions with its 
standard norm. Moreover, $M_\Phi$ is a Banach 
algebra with respect to pointwise operations and the norm
\[
\|a\|_{M_\Phi}
=
\|L(a)\|_{\mathcal{B}(\ell^\Phi(\mathbb{Z}))}
\]
(see \cite[Theorems~5.5--5.6 and Lemma~6.1]{KT-JAT}). 
This algebra contains all $2\pi$-periodic functions of bounded
variation (see \cite[Theorem~1.2(c)]{KT-JAT}).

Let $\mathbb{Z}_+:=\{0,1,2,\dots\}$ and let $P$ denote the discrete Riesz 
projection on $\ell^\Phi(\mathbb{Z})$ defined for 
$\varphi=\{\varphi_j\}_{j\in\mathbb{Z}}\in \ell^\Phi(\mathbb{Z})$
by
\[
(P\varphi)_j:=\left\{\begin{array}{ll}
0, & j\in\mathbb{Z}\setminus\mathbb{Z}_+,
\\
\varphi_j, &j\in\mathbb{Z}_+,
\end{array}\right.
\]
and $Q:=I-P$.
Consider the following subspaces of $\ell^\Phi(\mathbb{Z})$:
\begin{align*}
S_0(\mathbb{Z}_+) & :=PS_0(\mathbb{Z})=
\{\varphi=\{\varphi_j\}_{j\in\mathbb{Z}}\in S_0(\mathbb{Z}):\varphi_j=0
\mbox{ for }j\in\mathbb{Z}\setminus\mathbb{Z}_+ \},
\\
\ell^\Phi(\mathbb{Z}_+) &:=P\ell^\Phi(\mathbb{Z})=
\{\varphi=\{\varphi_j\}_{j\in\mathbb{Z}}\in \ell^\Phi(\mathbb{Z}):\varphi_j=0
\mbox{ for }j\in\mathbb{Z}\setminus\mathbb{Z}_+\}.
\end{align*}
Since the subspace $S_0(\mathbb{Z}_+)$ is dense in $\ell^\Phi(\mathbb{Z}_ +)$
in view of \cite[Lemma~3]{KT24}, for every $a\in M_\Phi$,
the operator from $S_0(\mathbb{Z}_+)$ to $\ell^\Phi(\mathbb{Z}_+)$ defined by 
$\varphi\mapsto P(a*\varphi)$, extends to a bounded operator
\[
T(a):\ell^\Phi(\mathbb{Z}_+)\to \ell^\Phi(\mathbb{Z}_+),
\quad
\varphi\mapsto P(a*\varphi),
\]
which is referred to as the discrete Wiener-Hopf (or Toeplitz) operator
with symbol $a$. It follows from \cite[Theorem~2]{KT24} that
if $a\in M_\Phi$, then
$\|T(a)\|_{\mathcal{B}(\ell^\Phi(\mathbb{Z}_+))}
= 
\|L(a)\|_{\mathcal{B}(\ell^\Phi(\mathbb{Z}))}$.

The Fredholm theory of discrete Wiener-Hopf operators on spaces
$\ell^p(\mathbb{Z}_+)$ and their weighted analogues
is a well developed branch of operator theory
(see, e.g., \cite[part~VI]{GGK93}, \cite[Ch.~III]{GGK03} for the case $p=2$ 
and 
\cite[Ch.~1]{BG05}, 
\cite[Ch.~7]{BS99}, 
\cite[Ch.~2 and~6]{BS06},
\cite{BSey00} 
for the case $1<p<\infty$). 
Gohberg and Krupnik proposed in 1973 a local principle, which  
became a convenient tool for studying Fredholm properties of 
convolution type operators. For expositions of the Gohberg-Krupnik
local principle, we refer to the English translation of the 
original 1973 paper \cite{GK10} and to 
\cite[Section~5.1]{GK92-vol-I}, \cite[Sections~1.30--1.32]{BS06},
\cite[Section~2.4]{RSS11}. Applications of the Gohberg-Krupnik
local principle to discrete Wiener-Hopf operators on $\ell^p$-spaces
were considered in \cite{GK10}, 
\cite[Section~9.14]{GK92-vol-II}, \cite[Theorem~2.69]{BS06},
\cite[Section~6]{BSey00}. The aim of this paper is to derive 
consequences of the Gohberg-Krupnik local principle for
discrete Wiener-Hopf operators in the setting of Orlicz sequence 
spaces and to extend \cite[Theorem~2.69]{BS06} to this setting.
In order to state our main result, we will need some notation.

Let $\Phi$ be an $N$-function with the Matuszewska-Orlicz indices 
satisfying
\begin{equation}\label{eq:indices-nontriviality}
1<\alpha_\Phi\le\beta_\Phi<\infty. 
\end{equation}
We refer to Section~\ref{sec:Matuszewska-Orlicz} and to 
\cite{MO60}, \cite[Section~3]{M85}, \cite[Ch.~11]{M89} for their 
definitions and properties. In particular, inequality 
\eqref{eq:indices-nontriviality} implies that the Orlicz space
$\ell^\Phi(\mathbb{Z})$ is {
reflexive and separable (reflexivity follows from inequalities
\eqref{eq:lower-upper-MO-indices} and 
Lemma~\ref{le:reflexivity-indices} below, separability follows
from reflexivity in view of 
\cite[Ch.~1, Corollaries 4.4 and~5.6]{BS88}).
Our results presented below, for instance, can be applied to 
Orlicz sequence spaces generated by $N$-functions
\[
\Phi(t)=t^p\ln^q (1+t),
\quad 
1<p<\infty, 
\quad
q\ge 0.
\]
It follows from \cite[Section~2, Example~7]{M85} that for 
such functions $\alpha_\Phi=p$ and $\beta_\Phi=p+q$.
}

For $\theta$ satisfying
\begin{equation}\label{eq:small-theta}
0<\theta<2\min\{1/\beta_\Phi,1-1/\alpha_\Phi\},
\end{equation}
let
\begin{equation}\label{eq:Phi-theta}
\Phi_\theta(t):=\widetilde{\Psi}_\theta^{-1}(t),
\quad
t\ge 0,
\end{equation}
where $\widetilde{\Psi}_\theta$ is the least
concave majorant of the function 
\begin{equation}\label{eq:Psi-theta}
\Psi_\theta(t):=
\left\{\begin{array}{cc}
\left(\Phi^{-1}(t)t^{-\frac{\theta}{2}}\right)^{\frac{1}{1-\theta}},
&
t> 0,
\\
0, & t=0.
\end{array}\right.
\end{equation}
One can show that $\Phi_\theta$
is an $N$-function whose Matuszewska-Orlicz indices satisfy
$1<\alpha_{\Phi_\theta}\le\beta_{\Phi_\theta}<\infty$
(see Theorem~\ref{th:interpolation-i} below). Let
\[
M_{[\Phi]}:=
\bigcap_{0<\theta<2\min\{1/\beta_\Phi,1-1/\alpha_\Phi\}}
M_{\Phi_\theta},
\quad
M_{\langle\Phi\rangle}:=
\bigcup_{0<\theta<2\min\{1/\beta_\Phi,1-1/\alpha_\Phi\}}
M_{\Phi_\theta}.
\]
Clearly, $M_{[\Phi]}\subset M_{\langle\Phi\rangle}$.
It follows from Lemma~\ref{le:inclusion} below that
$M_{\langle\Phi\rangle}\subset M_\Phi$.

Given $a\in L^\infty(-\pi,\pi)$ and an open finite interval
$u\subset\mathbb{R}$, we denote by $a|_u$ the restriction of 
$a$ to $u$ regarded as an element of $L^\infty(u)$. For 
$\tau\in[-\pi,\pi)$, let $\mathcal{U}_\tau$ denote
the family of all open finite intervals containing the point $\tau$.
The local distance between $a,b\in L^\infty(-\pi,\pi)$ at $\tau\in[-\pi,\pi)$
is defined by
\[
\operatorname{dist}_\tau(a,b):=
\inf_{u\in\mathcal{U}_\tau}\left\|a|_u-b|_u\right\|_{L^\infty(u)}.
\]
The following is our main result. It extends \cite[Theorem~2.69]{BS06}
from the setting of $\ell^p(\mathbb{Z}_+)$ to the setting of Orlicz
sequence spaces $\ell^\Phi(\mathbb{Z}_+)$ with the Matuszewska-Orlicz indices 
satisfying \eqref{eq:indices-nontriviality}.
\begin{theorem}\label{th:main}
Let $\Phi$ be an $N$-function with the Matuszewska-Orlicz indices
satisfying \eqref{eq:indices-nontriviality}.
Suppose that 
$a\in M_{\langle\Phi\rangle}$ and for each $\tau\in[-\pi,\pi)$, there
exists $a_\tau\in M_{[\Phi]}$ such that 
$\operatorname{dist}_\tau(a,a_\tau)=0$ and the discrete Wiener-Hopf
operator $T(a_\tau)$ is Fredholm on $\ell^\Phi(\mathbb{Z}_+)$. Then
the discrete Wiener-Hopf operator $T(a)$ is Fredholm on 
$\ell^\Phi(\mathbb{Z}_+)$.
\end{theorem}

{
Note that if $\Phi(t)=t^2$, then $\Phi_\theta(t)=t^2$ for all
$0<\theta<1$ and $M_{[\Phi]}=M_{\langle\Phi\rangle}=L^\infty(-\pi,\pi)$.
In this case Theorem~\ref{th:main} was obtained by Douglas and Sarason
\cite{DS70} by different methods.
}

The paper is organised as follows. In Section~\ref{sec:LWH}, we prove that
$\|L(a)\|_{\mathcal{B}(X(\mathbb{Z}))}=
\|L^\pi(a)\|_{\mathcal{B}(X(\mathbb{Z}))/\mathcal{K}(X(\mathbb{Z}))}$
and
$\|T(a)\|_{\mathcal{B}(X(\mathbb{Z}_+))}=
\|T^\pi(a)\|_{\mathcal{B}(X(\mathbb{Z}_+))/\mathcal{K}(X(\mathbb{Z}_+))}$
in the setting of reflexive translation-invariant Banach sequence
spaces $X(\mathbb{Z})$ and its subspaces $X(\mathbb{Z}_+):=PX(\mathbb{Z})$,
respectively.

In Section~\ref{sec:interpolation}, we collect definitions and some basic
properties of dilation exponents $\gamma_\Psi,\delta_\Psi$ of a function
$\Psi:(0,\infty)\to(0,\infty)$, as well as, of Matuszewska-Orlicz indices
of a $\varphi$-function $\Phi:[0,\infty)\to[0,\infty)$. Further, we show
if the Matuszewska-Orlicz indices of an $N$-function $\Phi$ satisfy
\eqref{eq:indices-nontriviality} and $\theta$ satisfies \eqref{eq:small-theta},
then $\Phi_\theta$, defined by \eqref{eq:Phi-theta}--\eqref{eq:Psi-theta},
is an $N$-function, whose Matsuszewska-Orlicz indices satisfy 
$1<\alpha_{\Phi_\theta}\le\beta_{\Phi_\theta}<\infty$, and a linear operator
$T$ bounded on $\ell^{\Phi_\theta}(\mathbb{Z})$ and on $\ell^2(\mathbb{Z})$
is also bounded on $\ell^\Phi(\mathbb{Z})$ with the norm estimate similar to 
that in the Riesz-Thorin interpolation theorem:
\begin{equation}\label{eq:interpolation-inequality}
\|T\|_{\mathcal{B}(\ell^\Phi(\mathbb{Z}))}
\le 
C_{\Phi,\theta}
\|T\|_{\mathcal{B}(\ell^{\Phi_\theta}(\mathbb{Z}))}^{1-\theta}
\|T\|_{\mathcal{B}(\ell^2(\mathbb{Z}))}^\theta.
\end{equation}
This is done by showing that the Orlicz sequence space $\ell^\Phi(\mathbb{Z})$
can be factorised as the Calder\'on product 
$(\ell^{\Phi_\theta}(\mathbb{Z}))^{1-\theta} (\ell^2(\mathbb{Z}))^\theta$,
and then applying an abstract interpolation theorem for Calder\'on products
(see \cite{C64} and \cite[Theorem~3.11]{MS19}).

In Section~\ref{sec:proof-main}, we state the abstract version of the 
Gohberg-Krupnik local principle \cite{GK10} and then show that it can
be applied in order to prove Theorem~\ref{th:main} following the scheme of 
the proof of \cite[Theorem~2.69]{BS06} (see also \cite[Section~6]{BSey00}).
The key point in the proof is to show that if 
$\operatorname{dist}_\tau(a,a_\tau)=0$, then the coset $T^\pi(a)$ is
$\mathcal{M}_\tau$-equivalent to its local representative $T^\pi(a_\tau)$, where $\mathcal{M}_\tau:=\{T^\pi(f):f\in\mathcal{N}_\tau\}$.
This is done by applying inequality \eqref{eq:interpolation-inequality}
to $L((a-a_\tau)f)$ for a family of ``bump" functions $f\in\mathcal{N}_\tau$.
The family $\mathcal{N}_\tau$ is chosen so that
$\sup_{f\in\mathcal{N}_\tau}
\|L((a-a_\tau)f)\|_{\mathcal{B}(\ell^{\Phi_\theta}(\mathbb{Z}))}
\le\operatorname{const}\|a-a_\tau\|_{M_{\Phi_{\theta}}}<\infty$ 
and
$\inf_{f\in\mathcal{N}_\tau}
\|L((a-a_\tau)f)\|_{\mathcal{B}(\ell^2(\mathbb{Z}))}
=\operatorname{dist}_\tau(a.a_\tau)=0$ for every 
$\tau\in[-\pi,\pi)$.
\section{Laurent and discrete Wiener-Hopf operators on translation-invariant
Banach sequence spaces}\label{sec:LWH}
\subsection{Banach sequence spaces and their associate spaces}
Let $\ell_+^0(\mathbb{Z})$ be the cone of nonnegative sequences in 
$\ell^0(\mathbb{Z})$. According to 
\cite[Ch.~1, Definition~1.1]{BS88}, a Banach function norm 
$\varrho:\ell_+^0(\mathbb{Z})\to [0,\infty]$ is a mapping which satisfies 
the following axioms for all $f,g\in \ell_+^0(\mathbb{Z})$, for all 
sequences $\{f^{(n)}\}_{n\in\mathbb{N}}$ in $\ell_+^0(\mathbb{Z})$, for all 
finite subsets $E\subset\mathbb{Z}$, and all constants $\alpha\ge 0$:
\begin{eqnarray*}
{\rm (A1)} & &
\varrho(f)=0  \Leftrightarrow  f=0,
\
\varrho(\alpha f)=\alpha\varrho(f),
\
\varrho(f+g) \le \varrho(f)+\varrho(g),\\
{\rm (A2)} & &0\le g \le f \ \ \Rightarrow \ 
\varrho(g) \le \varrho(f)
\quad\mbox{(the lattice property)},\\
{\rm (A3)} & &0\le f^{(n)} \uparrow f \ \ \Rightarrow \
       \varrho(f^{(n)}) \uparrow \varrho(f)\quad\mbox{(the Fatou property)},\\
{\rm (A4)} & & \varrho(\chi_E) <\infty,\\
{\rm (A5)} & &\sum_{k\in E} f_k \le C_E\varrho(f) ,
\end{eqnarray*}
where $\chi_E$ is the characteristic (indicator) function of $E$, and the 
constant $C_E \in (0,\infty)$ may depend on $\varrho$ and $E$, 
but is independent of $f$. The set $X(\mathbb{Z})$ of all sequences 
$f\in \ell^0(\mathbb{Z})$ for which $\varrho(|f|)<\infty$ is called a 
Banach sequence space. For each $f\in X(\mathbb{Z})$, the norm of $f$ is 
defined by 
\[
\|f\|_{X(\mathbb{Z})} :=\varrho(|f|). 
\]
The set $X(\mathbb{Z})$ equipped with the natural linear space operations 
and this norm becomes a Banach space 
(see \cite[Ch.~1, Theorems~1.4 and~1.6]{BS88}). If $\varrho$ 
is a Banach function norm, its associate norm $\varrho'$ is defined on 
$\ell_+^0(\mathbb{Z})$ by
\[
\varrho'(g):=\sup\left\{
\sum_{k\in\mathbb{Z}} f_kg_k \ : \ 
f=\{f_k\}_{k\in\mathbb{Z}}\in\ell_+^0(\mathbb{Z}), \ \varrho(f) \le 1
\right\}, \quad g\in \ell_+^0(\mathbb{Z}).
\]
It is a Banach function norm itself \cite[Ch.~1, Theorem~2.2]{BS88}.
The Banach sequence space $X'(\mathbb{Z})$ determined by the Banach 
function norm $\varrho'$ is called the associate space (K\"othe dual) 
of $X(\mathbb{Z})$. 
\subsection{The dual space of the space $X(\mathbb{Z})$ and 
its subspace $X(\mathbb{Z}_+)$} 
The associate space $X'(\mathbb{Z})$ can be viewed 
as a subspace of the Banach dual space $(X(\mathbb{Z}))^*$.
Consider the following subspaces
\begin{align*}
X(\mathbb{Z}_+) &:=PX(\mathbb{Z})=
\{\varphi=\{\varphi_j\}_{j\in\mathbb{Z}}\in X(\mathbb{Z}):\varphi_j=0
\mbox{ for }j\in\mathbb{Z}\setminus\mathbb{Z}_+\},
\\
X'(\mathbb{Z}_+) &:=PX'(\mathbb{Z})=
\{\varphi=\{\varphi_j\}_{j\in\mathbb{Z}}\in X'(\mathbb{Z}):\varphi_j=0
\mbox{ for }j\in\mathbb{Z}\setminus\mathbb{Z}_+\}
\end{align*}
of the space $X(\mathbb{Z})$ and its associate space $X'(\mathbb{Z})$,
respectively.

Now we are going to describe the dual space of a separable Banach 
sequence space $X(\mathbb{Z})$ and its subspace $X(\mathbb{Z}_+)$.
\begin{lemma}\label{le:dual-associate}
Let $X(\mathbb{Z})$ be a separable Banach sequence space. Then 
$(X(\mathbb{Z}))^*$ is isometrically isomorphic to $X'(\mathbb{Z})$ and
$(X(\mathbb{Z}_+))^*$ is isometrically isomorphic to $X'(\mathbb{Z}_+)$.
\end{lemma}
\begin{proof}
The fact that $(X(\mathbb{Z}))^*$ and $X'(\mathbb{Z})$ are isometrically
isomorphic follows from \cite[Ch.~1, Corollaries~4.3--4.4 and~5.6]{BS88}.
More precisely, for every $F\in (X(\mathbb{Z}))^*$ there is a unique
$y=\{y_n\}_{n\in\mathbb{Z}}\in X'(\mathbb{Z})$ such that for all
$x=\{x_n\}_{n\in\mathbb{Z}}\in X(\mathbb{Z})$,
\begin{equation}\label{eq:dual-associate-1}
F(x)=(x,y):=\sum_{n\in\mathbb{Z}} x_n\overline{y_n}
\end{equation}
and $\|F\|_{(X(\mathbb{Z}))^*}=\|y\|_{X'(\mathbb{Z})}$.

Let
\[
(X(\mathbb{Z}_+))^\perp:=\{F\in (X(\mathbb{Z}))^*\ :\ F(x)=0
\quad\mbox{for all}\quad x\in X(\mathbb{Z}_+)\}. 
\]
It follows from \eqref{eq:dual-associate-1} that
\[
(X(\mathbb{Z}_+))^\perp =(PX(\mathbb{Z}))^\perp= QX'(\mathbb{Z}).
\]
By \cite[Theorem~7.1]{D70}, $(X(\mathbb{Z}_+))^*$ is isometrically
isomorphic to 
$(X(\mathbb{Z}))^*/(X(\mathbb{Z}_+))^\perp=
X'(\mathbb{Z})/QX'(\mathbb{Z})$.
It remains to show that the latter space is isometrically isomorphic to
$X'(\mathbb{Z}_+)=PX'(\mathbb{Z})$.

It is easy to check that the mapping 
$\Psi:X'(\mathbb{Z})/QX'(\mathbb{Z})\to PX'(\mathbb{Z})$ given by
\begin{equation}\label{eq:dual-associate-2}
\Psi(x+QX'(\mathbb{Z}))=Px,
\quad
x\in X'(\mathbb{Z}),
\end{equation}
is an isomorphism. Moreover, if $x\in X'(\mathbb{Z})$, then
\begin{align}
\|x+QX'(\mathbb{Z})\|_{X'(\mathbb{Z})/QX'(\mathbb{Z})}
&=
\inf_{y\in QX'(\mathbb{Z})}\|x+y\|_{X'(\mathbb{Z})}
\nonumber\\
&=
\inf_{y\in QX'(\mathbb{Z})}\|Px+y\|_{X'(\mathbb{Z})}
\le 
\|Px\|_{X'(\mathbb{Z})}.
\label{eq:dual-associate-3}
\end{align}
On the other hand, for every 
$y=\{y_n\}_{n\in\mathbb{Z}}\in QX'(\mathbb{Z})$, the supports of 
$Px$ and $y$ are disjoint. Hence
\[
|(Px)_n|\le |(Px)_n+y_n|,
\quad
n\in\mathbb{Z},
\]
which implies that
\[
\|Px\|_{X'(\mathbb{Z})} \le \|Px+y\|_{X'(\mathbb{Z})},
\quad
y\in QX'(\mathbb{Z}).
\]
Therefore,
\begin{align}
\|Px\|_{X'(\mathbb{Z})} 
&\le 
\inf_{y\in QX'(\mathbb{Z})}\|Px+y\|_{X'(\mathbb{Z})}
\nonumber\\
&=
\inf_{y\in QX'(\mathbb{Z})}\|x+y\|_{X'(\mathbb{Z})}
=
\|x+QX'(\mathbb{Z})\|_{X'(\mathbb{Z})/QX'(\mathbb{Z})}.
\label{eq:dual-associate-4}
\end{align}
Combining \eqref{eq:dual-associate-3}--\eqref{eq:dual-associate-4},
we see that the isomorphism $\Psi$ given by 
\eqref{eq:dual-associate-2} is actually an isometry.
\end{proof}
\subsection{Strong convergence of shifted compact operators to the 
zero operator}
Let $X(\mathbb{Z})$ be a Banach sequence space. By $M_{X(\mathbb{Z})}$ 
we denote the collection of all distributions $a\in\mathcal{P}'$ for which 
$a*\varphi\in X(\mathbb{Z})$ whenever $\varphi\in S_0(\mathbb{Z})$ 
and
\[
\|a\|_{M_{X(\mathbb{Z})}}:=\sup\left\{
\frac{\|a*\varphi\|_{X(\mathbb{Z})}}
{\|\varphi\|_{X(\mathbb{Z})}}
\ :\
\varphi\in S_0(\mathbb{Z}), \ \varphi\ne 0
\right\}<\infty.
\]
If $X(\mathbb{Z})$ is separable, then $S_0(\mathbb{Z})$ is dense
in $X(\mathbb{Z})$ and {$S_0(\mathbb{Z}_+)$} is dense in $X(\mathbb{Z}_+)$
(see \cite[Lemma~3]{KT24}). So, for $a\in M_{X(\mathbb{Z})}$
the Laurent operator $L(a):S_0(\mathbb{Z})\to X(\mathbb{Z})$
defined by $L(a)\varphi=a*\varphi$ for $\varphi\in S_0(\mathbb{Z})$,
extends to a bounded operator
\[
L(a):X(\mathbb{Z})\to X(\mathbb{Z}),
\quad
\varphi\mapsto a*\varphi.
\]
Similarly, the discrete Wiener-Hopf operator 
$T(a):S_0(\mathbb{Z}_+) \to X(\mathbb{Z}_+)$ defined by 
$\varphi\mapsto P(a*\varphi)$, extends to a bounded operator
\[
T(a):X(\mathbb{Z}_+)\to X(\mathbb{Z}_+),
\quad
\varphi\mapsto P(a*\varphi).
\]

In order to guarantee that $M_{X(\mathbb{Z})}$ contains
nontrivial elements, we will impose an additional assumption
on a Banach sequence space $X(\mathbb{Z})$.
For $\varphi=\{\varphi_j\}_{j\in\mathbb{Z}}\in\ell^0(\mathbb{Z})$, define 
the translation (shift) operator 
\[
(T\varphi)_j=\varphi_{j-1},
\quad
j\in\mathbb{Z}.
\]
A Banach sequence space $X(\mathbb{Z})$ is said to be translation-invariant
if $T\varphi\in X(\mathbb{Z})$ and 
$\|T\varphi\|_{X(\mathbb{Z})}=\|\varphi\|_{X(\mathbb{Z})}$ for all
$\varphi=\{\varphi_j\}_{j\in\mathbb{Z}}\in X(\mathbb{Z})$.

If $X(\mathbb{Z})$ is translation-invariant, then
it is clear that $T^{-1}$ is a bounded operator on $X(\mathbb{Z})$ 
given by
\[
(T^{-1}\varphi)_j=\varphi_{j+1},
\quad
j\in\mathbb{Z}.
\]
As usual, $T^0:=I$ and $T^{-n}:=(T^{-1})^n$ for $n\in\mathbb{N}$.
It is easy to see that $\|T^m\|_{\mathcal{B}(X(\mathbb{Z}))}=1$
for all $m\in\mathbb{Z}$.

For $m\in\mathbb{Z}$, put
\[
\chi_m(\theta):=e^{im\theta},
\quad
\theta\in\mathbb{R}.
\]
If $\varphi\in S_0(\mathbb{Z})$, then
\[
(L(\chi_m)\varphi)_j
=
\sum_{k\in\mathbb{Z}}(\widehat{\chi_m})_{j-k}\varphi_k
=
\varphi_{j-m}
=(T^m\varphi)_j,
\quad
j\in\mathbb{Z}.
\]
So, if $m\in\mathbb{Z}$ and $X(\mathbb{Z})$ is a separable 
translation-invariant Banach sequence space, then 
$\chi_m\in M_{X(\mathbb{Z})}$ and
\[
\|\chi_m\|_{M_{X(\mathbb{Z})}}
=
\|L(\chi_m)\|_{\mathcal{B}(X(\mathbb{Z}))}
=
\|T^m\|_{\mathcal{B}(X(\mathbb{Z}))}=1.
\]
\begin{lemma}\label{le:shifted-compact-tend-to-zero}
Let $X(\mathbb{Z})$ be a reflexive translation-invariant Banach sequence space.
\begin{enumerate}
\item[{\rm(a)}]
If $K\in\mathcal{K}(X(\mathbb{Z}))$, then for all $\varphi\in X(\mathbb{Z})$,
\begin{equation}\label{eq:shifted-compact-tend-to-zero-1}
\lim_{n\to\infty}
\|L(\chi_{-n})KL(\chi_n)\varphi\|_{X(\mathbb{Z})}=0.
\end{equation}

\item[{\rm(b)}]
If $K\in\mathcal{K}(X(\mathbb{Z}_+))$, then for all 
$\varphi\in X(\mathbb{Z}_+)$,
\begin{equation}\label{eq:shifted-compact-tend-to-zero-2}
\lim_{n\to\infty}
\|T(\chi_{-n})KT(\chi_n)\varphi\|_{X(\mathbb{Z}_+)}=0.
\end{equation}
\end{enumerate}
\end{lemma}
\begin{proof}
It follows from Lemma~\ref{le:dual-associate} that $(X(\mathbb{Z}))^*$
is isometrically isomorphic to $X'(\mathbb{Z})$. In view of 
\cite[Ch.~1, Corollaries~4.3--4.4 and~5.6]{BS88}, the spaces $X(\mathbb{Z})$
and $X'(\mathbb{Z})$ are separable. Hence, by \cite[Lemma~3]{KT24},
$S_0(\mathbb{Z})$ is dense in $X(\mathbb{Z})$ and in $X'(\mathbb{Z})$,
and $S_0(\mathbb{Z})$ is dense in $X(\mathbb{Z}_+)$ and in $X'(\mathbb{Z}_+)$.

(a) Let 
$\varphi=\{\varphi_j\}_{j\in\mathbb{Z}}, \psi=\{\psi_j\}_{j\in\mathbb{Z}}
\in S_0(\mathbb{Z})$. If $n\in\mathbb{N}$ is large enough, then the supports of 
$L(\chi_n)\varphi$ and $\psi$ are disjoint, whence
\[
(L(\chi_n)\varphi,\psi)=
\sum_{j\in\mathbb{Z}}(L(\chi_n)\varphi)_j\overline{\psi_j}
=
\sum_{j\in\mathbb{Z}}\varphi_{j-n}\overline{\psi_j}=0.
\]
It follows from the above observation and \cite[Lemma~1.4.1(i)]{RSS11}
that $\{L(\chi_n)\}_{n\in\mathbb{Z}}$ converges weakly to the zero
operator on $X(\mathbb{Z})$. Since 
$\|L(\chi_{-n})\|_{\mathcal{B}(X(\mathbb{Z}))}=1$ for $n\in\mathbb{N}$,
the hypothesis $K\in\mathcal{K}(X(\mathbb{Z}))$ and 
\cite[Lemmas~1.4.4 and~1.4.6]{RSS11} imply that the sequence
$\{L(\chi_{-n})KL(\chi_n)\}_{n\in\mathbb{N}}$ converges strongly to the
zero operator on $X(\mathbb{Z})$, that is, 
\eqref{eq:shifted-compact-tend-to-zero-1} holds. Part(a) is proved.

(b) The proof is analogous to the proof of part (a). If 
$\varphi,\psi\in S_0(\mathbb{Z}_+)$ and $n\in\mathbb{N}$ is large enough, 
then the supports of $T(\chi_n)\varphi$ and $\psi$ are disjoint. 
Hence, in view of \cite[Lemma~1.4.1(i)]{RSS11}, the sequence tends
weakly to the zero operator on $X(\mathbb{Z}_+)$. Since 
$\|T(\chi_{-n})\|_{\mathcal{B}(X(\mathbb{Z}_+))}
\le\|L(\chi_{-n})\|_{\mathcal{B}(X(\mathbb{Z}))}=1$ for $n\in\mathbb{N}$,
as in the proof of part (a), we arrive at 
\eqref{eq:shifted-compact-tend-to-zero-2}.
\end{proof}
\subsection{Translation-invariance of Laurent and discrete Wiener-Hopf operators}
Let us show that $L(a)$ and $T(a)$ are invariant under translations
given by $L(\chi_{\pm n})$ and $T(\chi_{\pm n})$, respectively.
\begin{lemma}\label{le:TI-LWH}
Let $X(\mathbb{Z})$ be a separable translation-invariant Banach sequence space. Suppose $a\in M_{X(\mathbb{Z})}$. 
\begin{enumerate}
\item[{\rm(a)}]
For all $n\in\mathbb{Z}$, we have
\begin{equation}\label{eq:TI-LWH-1}
L(\chi_{-n})L(a)L(\chi_n)=L(a)
\end{equation}
on the space $X(\mathbb{Z})$.

\item[{\rm(b)}]
For all $n\in\mathbb{N}$, we have
\begin{equation}\label{eq:TI-LWH-2}
T(\chi_{-n})T(a)T(\chi_n)=T(a)
\end{equation}
on the space $X(\mathbb{Z}_+)$.
\end{enumerate}
\end{lemma}
\begin{proof}
(a) Let $n\in\mathbb{Z}$.
For $\varphi=\{\varphi_j\}_{j\in\mathbb{Z}}\in S_0(\mathbb{Z})$
and $j\in\mathbb{Z}$, we have
\begin{align*}
(L(\chi_{-n})L(a)L(\chi_n)\varphi)_j
&=
\sum_{k\in\mathbb{Z}}
\widehat{(\chi_{-n})}_{j-k}(L(a)L(\chi_n)\varphi)_k
=
(L(a)L(\chi_n)\varphi)_{j+n}
\\
&=
\sum_{k\in\mathbb{Z}}
\widehat{a}_{j+n-k}(L(\chi_n)\varphi)_k
=
\sum_{k\in\mathbb{Z}}
\widehat{a}_{j+n-k}
\left(
\sum_{m\in\mathbb{Z}}\widehat{(\chi_n)}_{k-m}\varphi_m
\right)
\\
&=
\sum_{k\in\mathbb{Z}}
\widehat{a}_{j+n-k}\varphi_{k-n}
=
\sum_{i\in\mathbb{Z}}\widehat{a}_{j-i}\varphi_i
=(L(a)\varphi)_j,
\end{align*}
that is, for all $\varphi\in S_0(\mathbb{Z})$,
\[
L(\chi_{-n})L(a)L(\chi_n)\varphi 
=
L(a)\varphi.
\]
Since $S_0(\mathbb{Z})$ is dense in $X(\mathbb{Z})$
(see \cite[Lemma~3]{KT24}), the above equality implies
\eqref{eq:TI-LWH-1} on the space $X(\mathbb{Z})$.
Part (a) is proved.

(b) Let $n\in\mathbb{N}$.
For $\varphi=\{\varphi_j\}_{j\in\mathbb{Z}}\in S_0(\mathbb{Z}_+)$
and $j\in\mathbb{Z}$, we have
\begin{align}
(T(\chi_{-n})T(a)T(\chi_n)\varphi)_j
&=
(PL(\chi_{-n})PL(a)PL(\chi_n)\varphi)_j
\nonumber\\
&=
\left\{\begin{array}{ll}
0, & j\in\mathbb{Z}\setminus\mathbb{Z}_+,
\\
(L(\chi_{-n})PL(a)PL(\chi_n)\varphi)_j, &j\in\mathbb{Z}_+
\end{array}\right.
\nonumber\\
&=
\left\{\begin{array}{ll}
0, & j\in\mathbb{Z}\setminus\mathbb{Z}_+,
\\
\sum_{k\in\mathbb{Z}}\widehat{(\chi_{-n})}_{j-k}
(PL(a)PL(\chi_n)\varphi)_k, & j\in\mathbb{Z}_+
\end{array}\right.
\nonumber\\
&=
\left\{\begin{array}{ll}
0, & j\in\mathbb{Z}\setminus\mathbb{Z}_+,
\\
(PL(a)PL(\chi_n)\varphi)_{j+n}, & j\in\mathbb{Z}_+.
\end{array}\right.
\label{eq:TI-LWH-3}
\end{align}
On the other hand, for all $j\in\mathbb{Z}$, we have
\begin{align*}
(PL(\chi_n)\varphi)_j
&=
\left\{\begin{array}{ll}
0, & j\in\mathbb{Z}\setminus\mathbb{Z}_+,
\\
(L(\chi_n)\varphi)_j, & j\in\mathbb{Z}_+
\end{array}\right.
=
\left\{\begin{array}{ll}
0, & j\in\mathbb{Z}\setminus\mathbb{Z}_+,
\\
\sum_{k\in\mathbb{Z}}\widehat{(\chi_n)}_{j-k}\varphi_k, 
& j\in\mathbb{Z}_+
\end{array}\right.
\\
&=
\left\{\begin{array}{ll}
0, & j\in\mathbb{Z}\setminus\mathbb{Z}_+,
\\
\varphi_{j-n}, 
& j\in\mathbb{Z}_+
\end{array}\right.
=
\left\{\begin{array}{ll}
0, & j-n\in\mathbb{Z}\setminus\mathbb{Z}_+,
\\
\varphi_{j-n}, 
& j-n\in\mathbb{Z}_+.
\end{array}\right.
\end{align*}
Hence, for $j\in\mathbb{Z}$, we obtain
\[
(L(a)PL(\chi_n)\varphi)_j
=
\sum_{k\in\mathbb{Z}}\widehat{a}_{j-k}(PL(\chi_n)\varphi)_k
=
\sum_{k=n}^\infty \widehat{a}_{j-k}\varphi_{k-n}
=
\sum_{m=0}^\infty\widehat{a}_{j-n-m}\varphi_m
\]
and
\[
(PL(a)PL(\chi_n)\varphi)_j
=
\left\{\begin{array}{ll}
0, & j\in\mathbb{Z}\setminus\mathbb{Z}_+,
\\
\sum_{k=0}^\infty\widehat{a}_{j-n-k}\varphi_k, &j\in\mathbb{Z}_+.
\end{array}\right.
\]
Thus
\begin{align}
(PL(a)PL(\chi_n)\varphi)_{j+n}
&=
\left\{\begin{array}{ll}
0, & j+n\in\mathbb{Z}\setminus\mathbb{Z}_+,
\\
\sum_{k=0}^\infty \widehat{a}_{j-k}\varphi_k, & j+n\in\mathbb{Z}_+
\end{array}\right.
\nonumber\\
&=
\left\{\begin{array}{ll}
0, & j+n\in\mathbb{Z}\setminus\mathbb{Z}_+,
\\
(L(a)\varphi)_j, & j+n\in\mathbb{Z}_+.
\end{array}\right.
\label{eq:TI-LWH-4}
\end{align}
Combining \eqref{eq:TI-LWH-3} and \eqref{eq:TI-LWH-4},
we see that for all $j\in\mathbb{Z}$,
\[
(T(\chi_{-n})T(a)T(\chi_n)\varphi)_j
=
\left\{\begin{array}{ll}
0, & j\in\mathbb{Z}\setminus\mathbb{Z}_+,
\\
(L(a)\varphi)_j, & j\in\mathbb{Z}_+
\end{array}\right.
=(PL(a)\varphi)_j=(T(a)\varphi)_j.
\]
Thus, for all $\varphi\in S_0(\mathbb{Z}_+)$, we get
\[
T(\chi_{-n})T(a)T(\chi_n)\varphi
=
T(a)\varphi.
\]
Since $S_0(\mathbb{Z}_+)$ is dense in $X(\mathbb{Z}_+)$
(see \cite[Lemma~3]{KT24}), the above equality yields
\eqref{eq:TI-LWH-2} on the space $X(\mathbb{Z}_+)$.
\end{proof}
\subsection{Maximal noncompactness of Laurent and discrete Wiener-Hopf operators}
Let $\mathcal{X}$ be a Banach space. An operator 
$A\in\mathcal{B}(\mathcal{X})$ is said to be maximally noncompact if 
$\|A\|_{\mathcal{B}(\mathcal{X})}=
\|A^\pi\|_{\mathcal{B}(\mathcal{X})/\mathcal{K}(\mathcal{X})}$.
Let us show that Laurent operators and discrete Wiener-Hopf operators
are maximally noncompact. 
\begin{theorem}
Let $X(\mathbb{Z})$ be a reflexive translation-invariant Banach sequence
space. If $a\in M_{X(\mathbb{Z})}$, then
\begin{align}
\|L^\pi(a)\|_{\mathcal{B}(X(\mathbb{Z}))/\mathcal{K}(X(\mathbb{Z}))}
&=
\|L(a)\|_{\mathcal{B}(X(\mathbb{Z}))},
\label{eq:norm=essential-norm-1}
\\
\|T^\pi(a)\|_{\mathcal{B}(X(\mathbb{Z}_+))/\mathcal{K}(X(\mathbb{Z}_+))}
&=
\|T(a)\|_{\mathcal{B}(X(\mathbb{Z}_+))}.
\label{eq:norm=essential-norm-2}
\end{align}
\end{theorem}
\begin{proof}
Equality \eqref{eq:norm=essential-norm-1} follows from Lemma~\ref{le:TI-LWH}(a),
\cite[Lemma~3]{KT24}, and \cite[Theorem~1.1]{KS24}. It can also be proved 
with the aid of Lemmas~\ref{le:shifted-compact-tend-to-zero}(a) 
and~\ref{le:TI-LWH}(a) by analogy with the proof of equality 
\eqref{eq:norm=essential-norm-2} presented below.

It follows from Lemma~\ref{le:TI-LWH}(b) and
\[
\|T(\chi_m)\|_{\mathcal{B}(X(\mathbb{Z}_+))}
\le
\|L(\chi_m)\|_{\mathcal{B}(X(\mathbb{Z}))}=1,
\quad
m\in\mathbb{Z},
\]
that for $\varphi\in X(\mathbb{Z}_+)$, $K\in\mathcal{K}(X(\mathbb{Z}_+))$,
and $n\in\mathbb{N}$, 
\begin{align*}
\|T(a)\varphi\|_{X(\mathbb{Z}_+)}
-
\|T(\chi_{-n})KT(\chi_n)\varphi\|_{X(\mathbb{Z}_+)}
&\le 
\|T(a)\varphi+T(\chi_{-n})KT(\chi_n)\varphi\|_{X(\mathbb{Z}_+)}
\\
&=
\|T(\chi_{-n})(T(a)+K)T(\chi_n)\varphi\|_{X(\mathbb{Z}_+)}
\\
&\le 
\|T(a)+K\|_{\mathcal{B}(X(\mathbb{Z}_+))}
\|\varphi\|_{X(\mathbb{Z}_+)}.
\end{align*}
Passing to the limit in the above inequality as $n\to\infty$,
we obtain in view of Lemma~\ref{le:shifted-compact-tend-to-zero}(b)
that for all $\varphi\in X(\mathbb{Z}_+)$ and 
$K\in\mathcal{K}(X(\mathbb{Z}_+))$,
\[
\|T(a)\varphi\|_{X(\mathbb{Z}_+)}
\le
\|T(a)+K\|_{\mathcal{B}(X(\mathbb{Z}_+))}
\|\varphi\|_{X(\mathbb{Z}_+)},
\]
whence
\begin{equation}\label{eq:norm=essential-norm-3}
\|T(a)\|_{\mathcal{B}(X(\mathbb{Z}_+))}
\le 
\inf_{K\in\mathcal{K}(X(\mathbb{Z}_+))}
\|T(a)+K\|_{\mathcal{B}(X(\mathbb{Z}_+))}
=
\|T^\pi(a)\|_{\mathcal{B}(X(\mathbb{Z}_+))/\mathcal{K}(X(\mathbb{Z}_+))}.
\end{equation}
The converse inequality is trivial:
\begin{equation}\label{eq:norm=essential-norm-4} 
\|T^\pi(a)\|_{\mathcal{B}(X(\mathbb{Z}_+))/\mathcal{K}(X(\mathbb{Z}_+))}
\le
\|T(a)\|_{\mathcal{B}(X(\mathbb{Z}_+))}.
\end{equation}
Combining inequalities 
\eqref{eq:norm=essential-norm-3}--\eqref{eq:norm=essential-norm-4},
we arrive at equality \eqref{eq:norm=essential-norm-2}.
\end{proof}
Since $L(a)$ and $T(a)$ are equal to zero if and only if their
symbol $a$ is trivial,
the above result immediately implies the following.
\begin{corollary}\label{co:noncompactness-Wiener-Hopf}
Let $X(\mathbb{Z})$ be a reflexive translation-invariant Banach sequence 
space. Suppose $a\in M_{X(\mathbb{Z})}$. Then 
\begin{enumerate}
\item[{\rm (a)}]
$L(a)\in\mathcal{K}(X(\mathbb{Z}))$ if and only if $a=0$;
\item[{\rm (b)}]
$T(a)\in\mathcal{K}(X(\mathbb{Z}_+))$ if and only if $a=0$.
\end{enumerate}
\end{corollary}

\subsection{Norms and essential norms of $L(a)$ and $T(a)$ are all equal
on a rearrangement-invariant Banach sequence space}
The distribution function of a 
sequence $f=\{f_k\}_{k\in\mathbb{Z}}\in\ell^0(\mathbb{Z})$ is defined by
\[
d_f(\lambda):=m\{k\in\mathbb{Z}:|f_k|>\lambda\},
\quad
\lambda\ge 0,
\]
where $m(S)$ denotes the measure (cardinality) of a set $S\subset\mathbb{Z}$.
One says that sequences $f=\{f_k\}_{k\in\mathbb{Z}}$, 
$g=\{g_k\}_{k\in\mathbb{Z}}\in\ell^0(\mathbb{Z})$ 
are equimeasurable if $d_f=d_g$. 
A Banach function norm $\varrho:\ell_+^0(\mathbb{Z})\to[0,\infty]$ is said to 
be rearrangement-invariant if $\varrho(f)=\varrho(g)$ for every pair of 
equimeasurable sequences 
$f=\{f_k\}_{k\in\mathbb{Z}}$, 
$g=\{g_k\}_{k\in\mathbb{Z}}\in\ell_+^0(\mathbb{Z})$. 
In that case, the Banach sequence space $X(\mathbb{Z})$ generated by $\varrho$ 
is said to be a rearrangement-invariant Banach sequence space 
(cf. \cite[Ch.~2, Definition~4.1]{BS88}).
It follows from \cite[Ch.~2, Proposition~4.2]{BS88} that if a Banach sequence 
space $X(\mathbb{Z})$ is rearrangement-invariant, then its associate space 
$X'(\mathbb{Z})$ is also a rearrangement-invariant Banach sequence space.
It is easy to see that for every 
$f=\{f_k\}_{k\in\mathbb{Z}}\in\ell^0_+(\mathbb{Z})$, the sequences $f$ 
and $Tf$ are equimeasurable. Hence, every rearrangement-invariant
Banach sequence space is also translation-invariant.
\begin{theorem}
Let $X(\mathbb{Z})$ be a reflexive rearrangement-invariant Banach sequence
space. If $a\in M_{X(\mathbb{Z})}$, then
\begin{align}
\|a\|_{M_{X(\mathbb{Z})}}
&=
\|L(a)\|_{\mathcal{B}(X(\mathbb{Z}))}
=
\|L^\pi(a)\|_{\mathcal{B}(X(\mathbb{Z}))/\mathcal{K}(X(\mathbb{Z}))}
\nonumber\\
&=
\|T(a)\|_{\mathcal{B}(X(\mathbb{Z}_+))}
=
\|T^\pi(a)\|_{\mathcal{B}(X(\mathbb{Z}_+))/\mathcal{K}(X(\mathbb{Z}_+))}.
\label{eq:corollary-norm=essential-norm}
\end{align}
\end{theorem}
\begin{proof}
The first equality in \eqref{eq:corollary-norm=essential-norm} is trivial. 
Since each rearrangement-invariant Banach sequence space is 
translation-invariant, the second and the fourth equalities in 
\eqref{eq:corollary-norm=essential-norm} are corollaries of 
equalities \eqref{eq:norm=essential-norm-1} and
\eqref{eq:norm=essential-norm-2}, respectively. Finally, the equality
$\|L(a)\|_{\mathcal{B}(X(\mathbb{Z}))}=
\|T(a)\|_{\mathcal{B}(X(\mathbb{Z}_+))}$ follows from a version of the 
Brown-Halmos theorem proved in \cite[Theorem~2]{KT24}.
\end{proof}

\section{Interpolation tools}\label{sec:interpolation}
\subsection{Dilation exponents and Matuszewska-Orlicz indices}
\label{sec:Matuszewska-Orlicz}
A function $\upsilon:(0,\infty)\to(0,\infty]$ is said to be submultiplicative
if $\upsilon(xy)\le\upsilon(x)\upsilon(y)$ for all $x,y\in(0,\infty)$. 
For a submultiplicative function $\upsilon$, put
\[
\alpha(\upsilon):=\sup_{0<t<1}\frac{\ln\upsilon(t)}{\ln t},
\quad
\beta(\upsilon):=\inf_{t>1}\frac{\ln\upsilon(t)}{\ln t}.
\]
\begin{theorem}\label{th:KPS-Maligranda}
Let $\upsilon$ be a submultiplicative function.
\begin{enumerate}
\item[{\rm(a)}]
If $\upsilon$ is everywhere finite, then
\begin{equation}\label{eq:inf-sup-lim}
\alpha(\upsilon)=\lim_{t\to 0}\frac{\ln\upsilon(t)}{\ln t},
\quad
\beta(\upsilon)=\lim_{t\to \infty}\frac{\ln\upsilon(t)}{\ln t}
\end{equation}
and $-\infty<\alpha(\upsilon)\le\beta(\upsilon)<\infty$.

\item[{\rm(b)}]
If $\upsilon$ is non-decreasing (not necessarily everywhere finite)
and $\upsilon(1)=1$, then \eqref{eq:inf-sup-lim} holds and
$0\le\alpha(\upsilon)\le\beta(\upsilon)\le\infty$.
\end{enumerate}
\end{theorem}
Part (a) is proved in \cite[Ch.II, Theorem~1.3]{KPS82} (see also
\cite[Theorem~1.13]{BK97}). Part (b) is contained in 
\cite[Remark~2 to Theorem~11.3]{M89}. 

Following \cite[Ch.~II, formula (1.16)]{KPS82}, for a function
$\Psi:(0,\infty)\to(0,\infty)$, its dilation function is defined by
\[
M(t,\Psi):=\sup_{0<u<\infty}\frac{\Psi(tu)}{\Psi(u)},
\quad
t\in(0,\infty).
\]
It is easy to see that it is submultiplicative. If it is everywhere
finite, then in view of Theorem~\ref{th:KPS-Maligranda}(a), there
exist numbers
\begin{equation}\label{eq:dilation-exponents}
\gamma_\Psi:=\alpha(M(\cdot,\Psi)),
\quad
\delta_\Psi:=\beta(M(\cdot,\Psi)),
\end{equation}
satisfying $-\infty<\gamma_\Psi\le\delta_\Psi<\infty$, which are
called the lower and the upper dilation exponents of $\Psi$,
respectively.

A function $\Phi:[0,\infty)\to[0,\infty)$ is called a $\varphi$-function
if it is continuous, nondecreasing, $\Phi(0)=0$, $\Phi(t)>0$ for $t>0$,
and $\Phi(t)\to\infty$ as $t\to\infty$. 
{It is clear that if $\Phi$ is an $N$-function, then
both $\Phi$ and its inverse $\Phi^{-1}$ are $\varphi$-functions.}

For a $\varphi$-function
$\Phi$, consider its dilation function $M(\cdot,\Phi)$, as well as,
the functions $M_0(\cdot,\Phi)$ and $M_\infty(\cdot,\Phi)$
defined for $t\in(0,\infty)$ by
\[
M_0(t,\Phi):=\limsup_{u\to 0^+}\frac{\Phi(tu)}{\Phi(u)},
\quad
M_\infty(t,\Phi):=\limsup_{u\to \infty}\frac{\Phi(tu)}{\Phi(u)}.
\] 
It is easy to check that the functions $M(\cdot,\Phi)$, $M_0(\cdot,\Phi)$,
and $M_\infty(\cdot,\Phi)$ are nondecreasing, submutliplicative, and have
value $1$  at $1$. Then 
\begin{align}
\alpha_\Phi:=\alpha(M(\cdot,\Phi)),
\quad
\alpha_\Phi^0:=\alpha(M_0(\cdot,\Phi)),
\quad
\alpha_\Phi^\infty:=\alpha(M_\infty(\cdot,\Phi)),
\label{eq:lower-MO-indices}
\\
\beta_\Phi:=\beta(M(\cdot,\Phi)),
\quad
\beta_\Phi^0:=\beta(M_0(\cdot,\Phi)),
\quad
\beta_\Phi^\infty:=\beta(M_\infty(\cdot,\Phi)),
\label{eq:upper-MO-indices}
\end{align}
are well defined. The quantities defined in \eqref{eq:lower-MO-indices}
are called the lower Matuszewska-Orlicz indices of $\Phi$ and the quantities
defined in \eqref{eq:upper-MO-indices} are called the upper 
Matuszewska-Orlicz indices of $\Phi$
(see \cite{MO60}, \cite[Section~3]{M85}, \cite[Ch.~11]{M89}). Since
\[
\max\{M_0(t,\Phi),M_\infty(t,\Phi)\}\le M(t,\Phi),
\quad
t\in (0,\infty),
\]
it follows from Theorem~\ref{th:KPS-Maligranda}(b) that
\begin{equation}\label{eq:lower-upper-MO-indices}
0\le\alpha_\Phi\le\min\{\alpha_\Phi^0,\alpha_\Phi^\infty\}
\le 
\max\{\beta_\Phi^0,\beta_\Phi^\infty\}\le\beta_\Phi\le\infty.
\end{equation}

The following well-known lemma is a consequence of 
\cite[Theorems~9.3,~11.7 and Corollary~11.6]{M89}.
\begin{lemma}\label{le:reflexivity-indices}
Let $\Phi$ be an $N$-function. Then the Orlicz sequence space 
$\ell^\Phi(\mathbb{Z})$ is reflexive if and only if the 
Matuszewska-Orlicz indices $\alpha_\Phi^0, \beta_\Phi^0$ of $\Phi$ satisfy 
\[
1<\alpha_\Phi^0\le\beta_\Phi^0<\infty.
\]
\end{lemma}
\subsection{The least concave majorant}
For a function $\Psi:[0,\infty)\to[0,\infty)$ such that $\Psi(t)=0$ if and 
only if $t=0$, the least concave majorant is defined by
\[
\widetilde{\Psi}(t):=\inf\Phi(t),
\quad
t\ge 0,
\]
where the infimum is taken over all functions $\Phi:[0,\infty)\to[0,\infty)$
satisfying
\[
{\Phi(t)=0}
\Leftrightarrow t=0;
\quad
\Phi \mbox{ is concave on } [0,\infty);
\quad
\Psi(t)\le\Phi(t) 
\mbox{ for all } 
t\in[0,\infty).
\]
It is known that 
\[
\widetilde{\Psi}(t)=\sup\sum_{i=1}^n \lambda_i\Psi(t_i),
\quad
t\ge 0,
\]
where the supremum is taken over all $n\in\mathbb{N}$ and all 
$\lambda_i,t_i\ge 0$ satisfying
\[
\sum_{i=1}^n \lambda_i=1,
\quad
\sum_{i=1}^n \lambda_it_i{=t}
\]
(see \cite[Ch.~II, formula (1.5)]{KPS82}). The following important lemma
is contained in \cite[Ch.~2, Corollary~2 to Lemma~1.2]{KPS82}.
\begin{lemma}\label{le:equivalence-to-concave-majorant}
Let $\Psi:[0,\infty)\to[0,\infty)$ be such that $\Psi(t)=0$ if and only if
$t=0$, and its dilation function 
$M(\cdot,\Psi)$ is everywhere finite. If $0<\gamma_\Psi\le\delta_\Psi<1$,
then $\Psi$ is equivalent to its least concave majorant $\widetilde{\Psi}$,
that is, there exists constants $C_1,C_2\in(0,\infty)$ such that
\[
C_1\Psi(t)\le\widetilde{\Psi}(t)\le C_2\Psi(t),
\quad
t\in[0,\infty).
\]
\end{lemma}
\subsection{Properties of the function $\Phi_\theta$}
Let us show that the function $\Phi_\theta$ inherits the properties
of the function $\Phi$ under assumption \eqref{eq:small-theta}.
\begin{theorem}\label{th:interpolation-i}
Let $\Phi$ be an $N$-function with the Matuszewska-Orlicz indices satisfying 
\eqref{eq:indices-nontriviality}. For $\theta$ satisfying 
\eqref{eq:small-theta}, let the function $\Phi_\theta$ be defined by 
\eqref{eq:Phi-theta}--\eqref{eq:Psi-theta}. Then it is an $N$-function with
the Matuszewska-Orlicz indices satisfying
\begin{equation}\label{eq:interpolation-i-0}
1<\alpha_{\Phi_\theta}\le\beta_{\Phi_\theta}<\infty.
\end{equation}
\end{theorem}

\begin{proof}
It follows from \eqref{eq:small-theta} and \cite[Theorem~11.5]{M89} that
\begin{equation}\label{eq:interpolation-proof-i-1}
0<\frac{\theta}{2}<\alpha_{\Phi^{-1}}<1,
\quad
0<\frac{\theta}{2}<1-\beta_{\Phi^{-1}}<1.
\end{equation}
If $M(s,\Phi^{-1})=\infty$ for some $s\in(0,\infty)$, then it follows 
from \cite[Theorem~11.2]{M89} and the fact that $M(\cdot,\Phi^{-1})$ 
is nondecreasing that $M(t,\Phi^{-1})=\infty$ for all $t>1$, which implies 
that $\beta_{\Phi^{-1}}=\infty$, but this is impossible since 
$\beta_{\Phi^{-1}}\in(0,1)$. Thus $M(s,\Phi^{-1})<\infty$ for all 
$s\in(0,\infty)$. In this case Theorem~\ref{th:KPS-Maligranda} implies that
\begin{equation}\label{eq:interpolation-proof-i-2}
\alpha_{\Phi^{-1}}=\gamma_{\Phi^{-1}},
\quad
\beta_{\Phi^{-1}}=\delta_{\Phi^{-1}}.
\end{equation}

If $\alpha_\Phi\le 2$, then $1-1/\alpha_\Phi\le 1/2$, whence
\begin{equation}\label{eq:interpolation-proof-i-3}
0<\theta<2\min\{1/\beta_\Phi,1-1/\alpha_\Phi\}\le 1.
\end{equation}
If $\alpha_\Phi>2$, then $1/\beta_\Phi\le 1/\alpha_\Phi\le 1/2$, 
whence \eqref{eq:interpolation-proof-i-3} holds. It follows from
\eqref{eq:interpolation-proof-i-1}--\eqref{eq:interpolation-proof-i-3}
that
\begin{align}
0 &<
\frac{1}{1-\theta}\left(\alpha_{\Phi^{-1}}-\frac{\theta}{2}\right)
=
\frac{1}{1-\theta}\left(\gamma_{\Phi^{-1}}-\frac{\theta}{2}\right),
\label{eq:interpolation-proof-i-4}
\\
1 &>
\frac{1}{1-\theta}\left(\beta_{\Phi^{-1}}-\frac{\theta}{2}\right)
=
\frac{1}{1-\theta}\left(\delta_{\Phi^{-1}}-\frac{\theta}{2}\right).
\label{eq:interpolation-proof-i-5}
\end{align}
Since the function $f(x)=x^{1/(1-\theta)}$ is strictly increasing,
it follows from \eqref{eq:Psi-theta} that
\begin{align}
M(s,\Psi_\theta)
&= 
\sup_{0<s<\infty}
\frac{
\left(
\Phi^{-1}(st)\cdot (st)^{-\frac{\theta}{2}}
\right)^{\frac{1}{1-\theta}}
}{
\left(
\Phi^{-1}(s)\cdot s^{-\frac{\theta}{2}}
\right)^{\frac{1}{1-\theta}}
}
=
\left(
\sup_{0<s<\infty}
\frac{\Phi^{-1}(st)}{\Phi^{-1}(s)}
\right)^{\frac{1}{1-\theta}}
\left(t^{-\frac{\theta}{2}}\right)^{\frac{1}{1-\theta}}
\nonumber\\
&=
\left(
M(t,\Phi^{-1})
\right)^{\frac{1}{1-\theta}}
\left(t^{-\frac{\theta}{2}}\right)^{\frac{1}{1-\theta}}.
\label{eq:interpolation-proof-i-6}
\end{align}
In particular, the function $M(\cdot,\Psi_\theta)$ is submutiplicative and
everywhere finite because $M(\cdot,\Phi^{-1})$ is so. It follows from 
Theorem~\ref{th:KPS-Maligranda}(a) and \eqref{eq:interpolation-proof-i-6}
that
\begin{align}
\gamma_{\Psi_\theta}
&=
\lim_{t\to 0^+} \frac{\ln M(t,\Psi_\theta)}{\ln t}
=
\frac{1}{1-\theta}\left(
\lim_{t\to 0^+}\frac{\ln M(t,\Phi^{-1})}{\ln t}
-
\frac{\theta}{2}\lim_{t\to 0^+}\frac{\ln t}{\ln t}
\right)
\nonumber\\
&=
\frac{1}{1-\theta}\left(\gamma_{\Phi^{-1}}-\frac{\theta}{2}\right)
\label{eq:interpolation-proof-i-7}
\end{align}
and
\begin{align}
\delta_{\Psi_\theta}
&=
\lim_{t\to \infty} \frac{\ln M(t,\Psi_\theta)}{\ln t}
=
\frac{1}{1-\theta}\left(
\lim_{t\to \infty}\frac{\ln M(t,\Phi^{-1})}{\ln t}
-
\frac{\theta}{2}\lim_{t\to \infty}\frac{\ln t}{\ln t}
\right)
\nonumber\\
&=
\frac{1}{1-\theta}\left(\delta_{\Phi^{-1}}-\frac{\theta}{2}\right).
\label{eq:interpolation-proof-i-8}
\end{align}
Combining 
\eqref{eq:interpolation-proof-i-4}--\eqref{eq:interpolation-proof-i-5} 
and
\eqref{eq:interpolation-proof-i-7}--\eqref{eq:interpolation-proof-i-8},
we get
\begin{equation}\label{eq:interpolation-proof-i-9}
0<\gamma_{\Psi_\theta} \le\delta_{\Psi_\theta}<1.
\end{equation}
Therefore, in view of Lemma~\ref{le:equivalence-to-concave-majorant},
there exist constants $C_{\theta,1},C_{\theta,2}\in(0,\infty)$
such that
\begin{equation}\label{eq:interpolation-proof-i-10}
C_{\theta,1} \Psi_\theta(t) 
\leq
\widetilde{\Psi}_\theta(t) 
\leq 
C_{\theta,2} \Psi_\theta(t),
\quad
t\in(0,\infty),
\end{equation}
where $\widetilde{\Psi}_\theta$ is the least concave majorant of $\Psi_\theta$.
Then
\begin{equation}\label{eq:interpolation-proof-i-11}
\frac{C_{\theta,1}}{C_{\theta,2}} M(t,\Psi_\theta) 
\le 
M(t,\widetilde{\Psi}_\theta)
\le 
\frac{C_{\theta,2}}{C_{\theta,1}} M(t,\Psi_\theta),
\quad
t\in(0,\infty).
\end{equation}
So $M(\cdot,\widetilde{\Psi}_\theta)$ is everywhere finite 
submutliplicative function. Then Theorem~\ref{th:KPS-Maligranda}(a)
and \eqref{eq:interpolation-proof-i-11} imply that
\begin{equation}\label{eq:interpolation-proof-i-12}
\gamma_{\Psi_\theta}=\gamma_{\widetilde{\Psi}_\theta},
\quad
\delta_{\Psi_\theta}=\delta_{\widetilde{\Psi}_\theta}.
\end{equation}

Let us show that $\widetilde{\Psi}_\theta$ is a $\varphi$-function.
It follows from the construction of $\widetilde{\Psi}_\theta$ that
$\widetilde{\Psi}(0)=0$, $\widetilde{\Psi}_\theta(t)>0$ for $t>0$, and
$\widetilde{\Psi}_\theta$ is concave. Then $\widetilde{\Psi}_\theta$
is nondecreasing and continuous on $(0,\infty)$
(see, e.g., \cite[Ch.~II, \S 1.1]{KPS82}). It remains to prove that 
$\widetilde{\Psi}_\theta(t)\to 0$ as $t\to 0^+$ and
$\widetilde{\Psi}_\theta(t)\to\infty$ as $t\to\infty$.

Fix $\varepsilon_1\in(0,\alpha_{\Phi^{-1}}-\theta/2)$. Then it follows 
from \eqref{eq:interpolation-proof-i-4} and \eqref{eq:lower-upper-MO-indices}
that
\[
0<\frac{1}{1-\theta}\left(
\alpha_{\Phi^{-1}}^0-\varepsilon_1-\frac{\theta}{2}
\right),
\quad
0<\frac{1}{1-\theta}\left(
\alpha_{\Phi^{-1}}^\infty-\varepsilon_1-\frac{\theta}{2}
\right).
\]
By \cite[Theorem~11.13]{M89}, there exist $q_0,q_\infty$ such that
\begin{align}
0<
\frac{1}{1-\theta}\left(
\alpha_{\Phi^{-1}}^0-\varepsilon_1-\frac{\theta}{2}
\right)
<
\frac{1}{1-\theta}\left(q_0-\frac{\theta}{2}\right)
\le 
\frac{1}{1-\theta}\left(
\alpha_{\Phi^{-1}}^0-\frac{\theta}{2}
\right),
\label{eq:interpolation-proof-i-13}
\\
0<
\frac{1}{1-\theta}\left(
\alpha_{\Phi^{-1}}^\infty-\varepsilon_1-\frac{\theta}{2}
\right)
<
\frac{1}{1-\theta}\left(q_\infty-\frac{\theta}{2}\right)
\le 
\frac{1}{1-\theta}\left(
\alpha_{\Phi^{-1}}^\infty-\frac{\theta}{2}
\right)
\label{eq:interpolation-proof-i-14}
\end{align}
and
\[
0<C_{q_0}:=
\sup_{0 < u,t \le 1}
\frac{\Phi^{-1}(tu)}{\Phi^{-1}(u)t^{q_0}}
<\infty,
\quad
0<C_{q_\infty}:=
\sup_{1\le u,t<\infty}
\frac{\Phi^{-1}(u)t^{q_\infty}}{\Phi^{-1}(tu)}
<\infty.
\]
Therefore
\begin{align}
&
\Phi^{-1}(t)\le C_{q_0}\Phi^{-1}(1)t^{q_0},
\quad
t\in(0,1],
\label{eq:interpolation-proof-i-15}
\\
&
C_{q_\infty}^{-1}\Phi^{-1}(1) t^{q_\infty} 
\le 
\Phi^{-1}(t),
\quad
t\in[1,\infty).
\label{eq:interpolation-proof-i-16}
\end{align}
It follows from \eqref{eq:interpolation-proof-i-10},
\eqref{eq:Psi-theta}, and
\eqref{eq:interpolation-proof-i-15} that for $t\in(0,1]$,
\[
\widetilde{\Psi}_\theta(t) 
\le 
C_{\theta,2}\left(
C_{q_0}\Phi^{-1}(1) t^{q_0}t^{-\frac{\theta}{2}}
\right)^{\frac{1}{1-\theta}}
=
C_{\theta,2}\left(C_{q_0}\Phi^{-1}(1)\right)^{\frac{1}{1-\theta}}
t^{\frac{1}{1-\theta}\left(q_0-\frac{\theta}{2}\right)}.
\]
This inequality and \eqref{eq:interpolation-proof-i-13} imply that
$\widetilde{\Psi}_\theta(t)\to 0$ as $t\to 0^+$.
Analogously, it follows from \eqref{eq:interpolation-proof-i-10},
 \eqref{eq:Psi-theta}, and
\eqref{eq:interpolation-proof-i-16} that for $t\in[1,\infty)$,
\[
\widetilde{\Psi}_\theta(t) 
\ge 
C_{\theta,1}\left(C_{q_\infty}^{-1}\Phi^{-1}(1) 
t^{q_\infty} t^{-\frac{\theta}{2}}\right)^{\frac{1}{1-\theta}}
=
C_{\theta,1}\left(C_{q_\infty}^{-1}\Phi^{-1}(1)\right)^{\frac{1}{1-\theta}}
t^{\frac{1}{1-\theta}\left(q_\infty-\frac{\theta}{2}\right)}.
\]
This inequality and \eqref{eq:interpolation-proof-i-14} imply that
$\widetilde{\Psi}_\theta(t)\to\infty$ as $t\to\infty$. Then
$\widetilde{\Psi}_\theta$ is a concave $\varphi$-function. 
Hence Theorem~\ref{th:KPS-Maligranda} implies that
\begin{equation}\label{eq:interpolation-proof-i-17}
\alpha_{\widetilde{\Psi}_\theta}
=
\gamma_{\widetilde{\Psi}_\theta},
\quad
\beta_{\widetilde{\Psi}_\theta}
=
\delta_{\widetilde{\Psi}_\theta}.
\end{equation}
Combining \eqref{eq:interpolation-proof-i-9},
\eqref{eq:interpolation-proof-i-12}, and
\eqref{eq:interpolation-proof-i-17}, we arrive at
\[
0<\alpha_{\widetilde{\Psi}_\theta}\le\beta_{\widetilde{\Psi}_\theta}<1.
\]

Since the function $\Phi_\theta$ defined by \eqref{eq:Phi-theta} is the 
convex $\varphi$-function inverse to the concave $\varphi$-function
$\widetilde{\Psi}_\theta$, it follows from the above inequalities and
\cite[Theorem~11.5]{M89} that inequalities \eqref{eq:interpolation-i-0}
hold.

It remains to prove that the convex $\varphi$-function $\Phi_\theta$ is, 
moreover, an $N$-function. Let $\varepsilon_2\in(0,\alpha_{\Phi_\theta}-1)$.
Then it follows from \eqref{eq:lower-upper-MO-indices} that
\[
1<\alpha_{\Phi_\theta}^0-\varepsilon_2,
\quad
1<\alpha_{\Phi_\theta}^\infty-\varepsilon_2.
\]
Then, by \cite[Theorem~11.13]{M89}, there exist $p_0,p_\infty$
such that
\begin{align}
&
1<\alpha_{\Phi_\theta}^0-\varepsilon_2<p_0\le \alpha_{\Phi_\theta}^0,
\label{eq:interpolation-proof-i-18}
\\
&
1<\alpha_{\Phi_\theta}^\infty-\varepsilon_2<p_\infty\le \alpha_{\Phi_\theta}^\infty
\label{eq:interpolation-proof-i-19}
\end{align}
and
\[
0<C_{p_0}:=\sup_{0<u,t\le 1}
\frac{\Phi_\theta(tu)}{\Phi_\theta(u)t^{p_0}}
<\infty,
\quad
0<C_{p_\infty}:=
\sup_{1\le u,t<\infty}
\frac{\Phi_\theta(u)t^{p_\infty}}{\Phi_\theta(tu)}
<\infty.
\]
Therefore
\begin{align}
&
\Phi_\theta(t)\le C_{p_0}\Phi_\theta(1)t^{p_0},
\quad
t\in(0,1],
\label{eq:interpolation-proof-i-20}
\\
&
C_{p_\infty}^{-1}\Phi_\theta(1) t^{p_\infty} \le \Phi_\theta(t),
\quad
t\in[1,\infty). 
\label{eq:interpolation-proof-i-21}
\end{align}
Taking into account \eqref{eq:interpolation-proof-i-18}
and \eqref{eq:interpolation-proof-i-20}, we get
\[
0\le\lim_{t\to 0^+}\frac{\Phi_\theta(t)}{t}
\le 
C_{p_0}\Phi_\theta(1)\lim_{t\to 0^+}t^{p_0-1}=0.
\]
Analogously, taking into account \eqref{eq:interpolation-proof-i-19}
and \eqref{eq:interpolation-proof-i-21}, we obtain
\[
\lim_{t\to\infty}\frac{\Phi_\theta(t)}{t}
\ge 
C_{p_\infty}^{-1}\Phi_\theta(1)\lim_{t\to\infty}t^{p_\infty-1}=\infty.
\]
Thus, $\Phi_\theta$ is an $N$-function.
\end{proof}
\subsection{Interpolation theorem of Riesz-Thorin type}
Let $\Phi_0,\Phi_1$ be $N$-functions and let
$\ell^{\Phi_0}(\mathbb{Z}),\ell^{\Phi_1}(\mathbb{Z})$ be 
the corresponding Orlicz sequence spaces. Suppose
$0<\theta<1$. The Calder\'on product 
$(\ell^{\Phi_0}(\mathbb{Z}))^{1-\theta} 
(\ell^{\Phi_1}(\mathbb{Z}))^\theta$
of $\ell^{\Phi_0}(\mathbb{Z})$ and $\ell^{\Phi_1}(\mathbb{Z})$ consists of 
all sequences
$x=\{x_k\}_{k\in\mathbb{Z}}\in\ell^0(\mathbb{Z})$ such that the inequality
\begin{equation}\label{eq:Calderon-product}
|x_k|\le\lambda |y_k|^{1-\theta}|z_k|^\theta,
\quad
k\in\mathbb{Z},
\end{equation}
holds for some $\lambda>0$ and elements 
$y=\{y_k\}_{k\in\mathbb{Z}}\in \ell^{\Phi_0}(\mathbb{Z})$ and
$z=\{z_k\}_{k\in\mathbb{Z}}\in \ell^{\Phi_1}(\mathbb{Z})$ with
$\|y\|_{\ell^{\Phi_0}(\mathbb{Z})}\le 1$ and 
$\|z\|_{\ell^{\Phi_1}(\mathbb{Z})}\le 1$
(see \cite{C64}).
The norm of $x$ in 
$(\ell^{\Phi_0}(\mathbb{Z}))^{1-\theta} 
(\ell^{\Phi_1}(\mathbb{Z}))^\theta$
is defined to be the infimum of all values of $\lambda$ in
\eqref{eq:Calderon-product}.
\begin{theorem}\label{th:interpolation-ii}
Let $\Phi$ be an $N$-function with the Matuszewska-Orlicz indices satisfying 
\eqref{eq:indices-nontriviality}. For $\theta$ satisfying 
\eqref{eq:small-theta}, let the function $\Phi_\theta$ be defined by 
\eqref{eq:Phi-theta}--\eqref{eq:Psi-theta}. Then
there exists a constant $C_{\Phi,\theta}\in(0,\infty)$ such that
for every linear operator 
\[
T:\ell^{\Phi_\theta}(\mathbb{Z})+\ell^2(\mathbb{Z})\to
\ell^{\Phi_\theta}(\mathbb{Z})+\ell^2(\mathbb{Z})
\]
such that $T\in\mathcal{B}(\ell^{\Phi_\theta}(\mathbb{Z}))$ and 
$T\in\mathcal{B}(\ell^2(\mathbb{Z}))$, one has 
$T\in\mathcal{B}(\ell^\Phi(\mathbb{Z}))$ and
\begin{equation}\label{eq:interpolation-ii-0}
\|T\|_{\mathcal{B}(\ell^\Phi(\mathbb{Z}))}
\le 
C_{\Phi,\theta}
\|T\|_{\mathcal{B}(\ell^{\Phi_\theta}(\mathbb{Z}))}^{1-\theta}
\|T\|_{\mathcal{B}(\ell^2(\mathbb{Z}))}^\theta.
\end{equation}
\end{theorem}
\begin{proof}
The argument below is a minor modification of the argument in 
\cite[Ch.~14, Example~3, pp.~178--180]{M89}. For the convenience of the reader,
we provide a detailed proof. Let
\[
\Phi_0(t):=\Phi_\theta(t),
\quad
\Phi_1(t):=t^2,
\quad
\rho(s,t):=s^{1-\theta}t^\theta,
\quad
s,t\ge 0.
\]
It follows from \eqref{eq:Psi-theta} that
\[
\Phi^{-1}(t)=(\Psi_\theta(t))^{1-\theta}(t^{1/2})^\theta,
\quad
t\ge 0.
\]
This equality, definition \eqref{eq:Phi-theta}, and inequality 
\eqref{eq:interpolation-proof-i-10} imply that
\[
C_{\theta,1}^{1-\theta}\Phi^{-1}(t)
\le 
(\Phi_\theta^{-1}(t))^{1-\theta}(t^{1/2})^\theta
\le 
C_{\theta,2}^{1-\theta}\Phi^{-1}(t),
\quad
t\ge 0,
\]
that is,
\begin{equation}\label{eq:interpolation-proof-ii-1}
C_{\theta,1}^{1-\theta}\Phi^{-1}(t)
\le 
\rho(\Phi_0^{-1}(t),\Phi_1^{-1}(t))
\le 
C_{\theta,2}^{1-\theta}\Phi^{-1}(t),
\quad
t\ge 0.
\end{equation}
If
\[
x=\{x_j\}_{j\in\mathbb{Z}} \in X_\theta(\mathbb{Z})
:=
(\ell^{\Phi_0}(\mathbb{Z}))^{1-\theta}
(\ell^{\Phi_1}(\mathbb{Z}))^\theta
\]
and
$\|x\|_{X_\theta(\mathbb{Z})}<\lambda$, then there exist
$x^{(i)}=\{x_j^{(i)}\}_{j\in\mathbb{Z}}\in\ell^{\Phi_i}(\mathbb{Z})$,
$i=0,1$, such that
\[
\|x^{(i)}\|_{\ell^{\Phi_i}(\mathbb{Z})}\le 1,
\quad
i=0,1
\]
and
\begin{equation}\label{eq:interpolation-proof-ii-2}
|x_j|\le\lambda\,\rho\left(|x_j^{(0)}|,|x_j^{(1)}|\right),
\quad
j\in\mathbb{Z}.
\end{equation}
Put
\begin{equation}\label{eq:interpolation-proof-ii-3}
y_j^{(i)}:=\Phi_i(|x_j^{(i)}|),
\quad
i=0,1,
\quad
j\in\mathbb{Z}.
\end{equation}
\begin{equation}\label{eq:interpolation-proof-ii-4}
y_j:=\max\{y_j^{(0)},y_j^{(1)}\},
\quad
j\in\mathbb{Z}.
\end{equation}
It follows from 
\eqref{eq:interpolation-proof-ii-1}--\eqref{eq:interpolation-proof-ii-4}
that for all $j\in\mathbb{Z}$,
\begin{align*}
\Phi\left(\frac{|x_j|}{\lambda C_{\theta,2}^{1-\theta}}\right)
&\le 
\Phi\left(
\frac{1}{C_{\theta,2}^{1-\theta}} 
\rho\left(|x_j^{(0)}|,|x_j^{(1)}|\right)
\right)
= 
\Phi\left(
\frac{1}{C_{\theta,2}^{1-\theta}} 
\rho\left(\Phi_0^{-1}(y_j^{(0)}),\Phi_1^{-1}(y_j^{(1)})\right)
\right)
\\
&\le
\Phi\left(
\frac{1}{C_{\theta,2}^{1-\theta}} 
\rho\left(\Phi_0^{-1}(y_j),\Phi_1^{-1}(y_j)\right)
\right)
\le 
\Phi\left(\Phi^{-1}(y_j)\right)=y_j.
\end{align*}
Taking into account that $\Phi$ is convex and $\Phi(0)=0$,
we obtain from the above inequality and 
\eqref{eq:interpolation-proof-ii-3}--\eqref{eq:interpolation-proof-ii-4}
that
\begin{align*}
\varrho_\Phi\left(\frac{x}{2\lambda C_{\theta,2}^{1-\theta}}\right)
&\le 
\frac{1}{2}\varrho_\Phi\left(\frac{x}{\lambda C_{\theta,2}^{1-\theta}}\right)
\le 
\frac{1}{2}\sum_{j\in\mathbb{Z}}y_j 
\\
&\le 
\frac{1}{2}\sum_{j\in\mathbb{Z}}y_j^{(0)}
+
\frac{1}{2}\sum_{j\in\mathbb{Z}}y_j^{(1)}
=
\frac{1}{2}\varrho_{\Phi_0}(x^{(0)})
+
\frac{1}{2}\varrho_{\Phi_1}(x^{(1)}).
\end{align*}
Since $\|x^{(i)}\|_{\ell^{\Phi_i}(\mathbb{Z})}\le 1$, $i=0,1$, we conclude that
\[
\frac{1}{2}\varrho_{\Phi_0}(x^{(0)})
+
\frac{1}{2}\varrho_{\Phi_1}(x^{(1)})
\le 
\frac{1}{2}+\frac{1}{2}=1.
\]
Therefore 
\[
\varrho_\Phi\left(\frac{x}{2\lambda C_{\theta,2}^{1-\theta}}\right)\le 1,
\]
which implies that
\[
\left\|
\frac{x}{2\lambda C_{\theta,2}^{1-\theta}}
\right\|_{\ell^\Phi(\mathbb{Z})} \le 1.
\]
Hence $\|x\|_{\ell^\Phi(\mathbb{Z})} \le 2C_{\theta,2}^{1-\theta}\lambda$
for  all $\lambda>\|x\|_{X_\theta(\mathbb{Z})}$. Thus
\begin{equation}\label{eq:interpolation-proof-ii-5}
\|x\|_{\ell^\Phi(\mathbb{Z})}
\le 
2C_{\theta,2}^{1-\theta}
\inf\{\lambda:\lambda>\|x\|_{X_\theta(\mathbb{Z})}\}
=
2C_{\theta,2}^{1-\theta}\|x\|_{X_\theta(\mathbb{Z})}.
\end{equation}
On the other hand, if $x\in\ell^\Phi(\mathbb{Z})\setminus\{0\}$, then it 
follows from \eqref{eq:interpolation-proof-ii-1} that for all $j\in\mathbb{Z}$,
\begin{equation}\label{eq:interpolation-proof-ii-6}
\frac{C_{\theta,1}^{1-\theta}|x_j|}{\|x\|_{\ell^\Phi(\mathbb{Z})}}
=
C_{\theta,1}^{1-\theta}
\Phi^{-1}\left(\Phi\left(
\frac{|x_j|}{\|x\|_{\ell^\Phi(\mathbb{Z})}}
\right)\right)
\le
\rho\left(y_j^{(0)},y_j^{(1)}\right),
\end{equation}
where
\begin{equation}\label{eq:interpolation-proof-ii-7}
y_j^{(i)}:=\Phi_i^{-1}\left(\Phi\left(
\frac{|x_j|}{\|x\|_{\ell^\Phi(\mathbb{Z})}}
\right)\right),
\quad
i=0,1,
\quad
j\in\mathbb{Z}.
\end{equation}
Since $x\in\ell^\Phi(\mathbb{Z})$, one has
\[
\varrho_\Phi\left(\frac{x}{\|x\|_{\ell^\Phi(\mathbb{Z})}}\right) 
\le 1.
\]
Taking into account the above inequality and 
\eqref{eq:interpolation-proof-ii-7}, we obtain
\[
\varrho_{\Phi_i}(y^{(i)})
=
\varrho_\Phi\left(\frac{x}{\|x\|_{\ell^\Phi(\mathbb{Z})}}\right) 
\le 1,
\quad
i=0,1,
\]
which implies that $\|y^{(i)}\|_{\ell^{\Phi_i}(\mathbb{Z})}\le 1$
for $i=0,1$. It follows from \eqref{eq:interpolation-proof-ii-6} that
\[
|x_j|\le\frac{\|x\|_{\ell^\Phi(\mathbb{Z})}}{C_{\theta,1}^{1-\theta}}
\rho\left(y_j^{(0)},y_j^{(1)}\right),
\quad
j\in\mathbb{Z}.
\]
Hence $x\in X_\theta(\mathbb{Z})$ and
\begin{equation}\label{eq:interpolation-proof-ii-8}
\|x\|_{X_\theta(\mathbb{Z})} 
\le 
\frac{\|x\|_{\ell^\Phi(\mathbb{Z})}}{C_{\theta,1}^{1-\theta}}.
\end{equation}
Thus \eqref{eq:interpolation-proof-ii-5} and \eqref{eq:interpolation-proof-ii-8}
imply that $X_\theta(\mathbb{Z})=\ell^\Phi(\mathbb{Z})$ and
for all $x\in\ell^\Phi(\mathbb{Z})$,
\begin{equation}\label{eq:interpolation-proof-ii-9}
C_{\theta,1}^{1-\theta} \|x\|_{X_\theta(\mathbb{Z})}
\le 
\|x\|_{\ell^\Phi(\mathbb{Z})}
\le 
2C_{\theta,2}^{1-\theta} \|x\|_{X_\theta(\mathbb{Z})}.
\end{equation}

It follows from \cite[Theorem~3.11]{MS19} that if 
$T:\ell^{\Phi_0}(\mathbb{Z})+\ell^{\Phi_1}(\mathbb{Z})\to
\ell^{\Phi_0}(\mathbb{Z})+\ell^{\Phi_1}(\mathbb{Z})$ is such that
$T\in\mathcal{B}(\ell^{\Phi_i}(\mathbb{Z}))$, $i=0,1$, then
$T\in\mathcal{B}(X_\theta(\mathbb{Z}))$ and
\begin{equation}\label{eq:interpolation-proof-ii-10}
\|T\|_{\mathcal{B}(X_\theta(\mathbb{Z}))}
\le 
\|T\|_{\mathcal{B}(\ell^{\Phi_0}(\mathbb{Z}))}^{1-\theta}
\|T\|_{\mathcal{B}(\ell^{\Phi_1}(\mathbb{Z}))}^{\theta}.
\end{equation}
If $x\in\ell^\Phi(\mathbb{Z})$, then it follows from 
\eqref{eq:interpolation-proof-ii-9}--\eqref{eq:interpolation-proof-ii-10}
that
\begin{align*}
\|Tx\|_{\ell^\Phi(\mathbb{Z})}
&\le 
2C_{\theta,2}^{1-\theta}\|Tx\|_{X_\theta(\mathbb{Z})}
\le 
2C_{\theta,2}^{1-\theta}
\|T\|_{\mathcal{B}(\ell^{\Phi_0}(\mathbb{Z}))}^{1-\theta}
\|T\|_{\mathcal{B}(\ell^{\Phi_1}(\mathbb{Z}))}^{\theta}
\|x\|_{X_\theta(\mathbb{Z})}
\\
&\le
\frac{2C_{\theta,2}^{1-\theta}}{C_{\theta,1}^{1-\theta}}
\|T\|_{\mathcal{B}(\ell^{\Phi_0}(\mathbb{Z}))}^{1-\theta}
\|T\|_{\mathcal{B}(\ell^{\Phi_1}(\mathbb{Z}))}^{\theta}
\|x\|_{\ell^\Phi(\mathbb{Z})}.
\end{align*}
The above inequality implies that \eqref{eq:interpolation-ii-0}
holds with $C_{\Phi,\theta}:=2C_{\theta,2}^{1-\theta} /C_{\theta,1}^{1-\theta}$,
where $C_{\theta,1},C_{\theta,2}$ are the constants from 
\eqref{eq:interpolation-proof-i-10}.

\end{proof}
\section{Proof of the main result}\label{sec:proof-main}
\subsection{The Gohberg-Krupnik local principle}
Let $\mathcal{A}$ be a complex Banach algebra with identity $e$ and let 
$\mathcal{T}$ be an arbitrary set. Suppose for each $\tau\in\mathcal{T}$
we are given a set $\mathcal{M}_\tau\subset\mathcal{A}$ such that
$0\notin\mathcal{M}_\tau$ and such that for every $f_1,f_2\in\mathcal{M}_\tau$
there exists an element $f\in\mathcal{M}_\tau$ such that
\[
f_jf=ff_j=f,
\quad j=1,2.
\] 
Suppose that every collection $\{f_\tau\}_{\tau\in\mathcal{T}}$
of elements $f_\tau\in\mathcal{M}_\tau$ contains a finite subcollection
$\{f_{\tau_j}\}_{j=1}^k$ such that $\sum_{j=1}^k f_{\tau_j}$ is invertible
in $\mathcal{A}$. Put
\[
\mathcal{M}=\bigcup_{\tau\in\mathcal{T}}\mathcal{M}_\tau.
\]
Denote by $\operatorname{Com}\mathcal{M}$ the set of all $a\in\mathcal{A}$
such that
\[
am=ma\quad\mbox{for all}\quad m\in\mathcal{M}.
\]
Two elements $a,b\in\mathcal{A}$ are said to be $\mathcal{M}_\tau$-equivalent 
if
\[
\inf_{f\in\mathcal{M}_\tau}\|(a-b)f\|_{\mathcal{A}}
=
\inf_{f\in\mathcal{M}_\tau}\|f(a-b)\|_{\mathcal{A}}=0.
\]
\begin{theorem}[Gohberg and Krupnik]
\label{th:Gohberg-Krupnik}
If $a\in\operatorname{Com}\mathcal{M}$ and if for every $\tau\in\mathcal{T}$
there exists an invertible element $a_\tau\in\mathcal{A}$ such that
$a$ and $a_\tau$ are $\mathcal{M}_\tau$-equivalent, then $a$ itself is 
invertible.
\end{theorem}
We refer to \cite[Theorem~1.1]{GK10}, 
\cite[Ch.~5, Theorem~1.1]{GK92-vol-I}, \cite[Theorem-1.32]{BS06},
\cite[Theorem~2.4.5]{RSS11} for the proof of this result.
\subsection{Discrete Wiener-Hopf operators and Hankel operators}
Let $\Phi$ be an $N$-function. It is easy to see that the flip operator
$J$ defined for $\varphi=\{\varphi_n\}_{n\in\mathbb{Z}}\in\ell^\Phi(\mathbb{Z})$
by
\[
(J\varphi)_n=\varphi_{-n-1},
\quad
n\in\mathbb{Z},
\]
is an isometry on $\ell^\Phi(\mathbb{Z})$. Assume that $\ell^\Phi(\mathbb{Z})$
is separable and put $Q:=I-P$. For $a\in M_\Phi$, define the Hankel operators
\begin{align*}
&
H(a):\ell^\Phi(\mathbb{Z}_+)\to\ell^\Phi(\mathbb{Z}_+),
\quad
\varphi\mapsto PL(a)QJ\varphi,
\\
&
H(\widetilde{a}):\ell^\Phi(\mathbb{Z}_+)\to\ell^\Phi(\mathbb{Z}_+),
\quad
\varphi\mapsto JQL(a)P\varphi,
\end{align*}
generated by $a$. Both $H(a)$ and $H(\widetilde{a})$ are bounded on the
subspace $\ell^\Phi(\mathbb{Z}_+)$ and their action on $\ell^\Phi(\mathbb{Z}_+)$
is given by the formulas
\begin{align}
&
H(a):\ell^\Phi(\mathbb{Z}_+)\to\ell^\Phi(\mathbb{Z}_+),
\quad
\{\varphi_j\}_{j\in\mathbb{Z}_+}
\mapsto
\left\{
\sum_{k\in\mathbb{Z}_+}\widehat{a}_{j+k+1}\varphi_k
\right\}_{j\in\mathbb{Z}_+},
\label{eq:Ha-representation}
\\
&
H(\widetilde{a}):\ell^\Phi(\mathbb{Z}_+)\to\ell^\Phi(\mathbb{Z}_+),
\quad
\{\varphi_j\}_{j\in\mathbb{Z}_+}
\mapsto
\left\{
\sum_{k\in\mathbb{Z}_+}\widehat{a}_{-j-k-1}\varphi_k
\right\}_{j\in\mathbb{Z}_+}
\label{eq:Ha-tilde-representation}
\end{align}
(recall that $\varphi_j=0$ if $j\in\mathbb{Z}\setminus\mathbb{Z}_+$).

The following lemma is an extension of \cite[Proposition~2.14]{BS06} from
the setting of $\ell^p$-spaces, $1<p<\infty$, to the setting of reflexive
Orlicz sequence spaces, with the same proof.
\begin{lemma}\label{le:Widom-formulas}
Let $\Phi$ be an $N$-function such that the corresponding Orlicz space
$\ell^\Phi(\mathbb{Z})$ is reflexive. If $a,b\in M_\Phi$, then
\[
T(ab)=T(a)T(b)+H(a)H(\widetilde{b})
\]
on the subspace $\ell^\Phi(\mathbb{Z}_+)$.
\end{lemma}
\subsection{The algebra $C_\Phi$ of continuous multipliers}
Let $\Phi$ be an $N$-function such that the corresponding Orlicz sequence
space $\ell^\Phi(\mathbb{Z})$ is reflexive.
Let $\mathcal{TP}$ denote the set of all trigonometric polynomials,
that is, functions of the form
\[
p(\theta)=\sum_{|k|\le n} c_k e^{ik\theta},
\quad\theta\in\mathbb{R},
\]
where $n\in\mathbb{Z}_+$ and $c_k\in\mathbb{C}$ for all $k\in\{-n,\dots,n\}$.
It follows from \cite[Lemma~3.1]{L68} (see also \cite[Theorem~6]{KT24}) that
$\mathcal{TP}\subset M_\Phi$. We denote by 
$C_\Phi$ the closure of $\mathcal{TP}$ with respect 
to the norm of the Banach algebra $M_\Phi$.

Let $C[-\pi,\pi]$ denote the set of all continuous $2\pi$-periodic
functions $f:\mathbb{R}\to\mathbb{C}$ with the supremum norm.
A $2\pi$-periodic function $a:\mathbb{R}\to\mathbb{C}$ is said to 
be of bounded variation if
\[
V(a):=\sup\sum_{j=1}^n |a(\theta_{j+1})-a(\theta_j)|<\infty,
\]
where the supremum is taken over all partitions 
$-\pi\le\theta_1<\theta_2<\dots<\theta_{n+1}=\pi$. It is well known
that the collection of all functions of bounded variation $BV[-\pi,\pi]$
is a Banach algebra with respect to the norm
\[
\|a\|_{BV[-\pi,\pi]}:=\|a\|_{L^\infty(-\pi,\pi)}+V(a).
\]
The following lemma is analogous to \cite[Lemma~6.2]{BSey00}.
\begin{lemma}\label{le:CBV-in-C-Phi}
Let $\Phi$ be an $N$-function such that the corresponding Orlicz sequence
space $\ell^\Phi(\mathbb{Z})$ is reflexive. Then
\[
C[-\pi,\pi]\cap BV[-\pi,\pi]\subset C_\Phi\subset C[-\pi,\pi].
\]
\end{lemma}
\begin{proof}
It follows from Lemma~\ref{le:reflexivity-indices}, 
\cite[Theorem~11.5]{M89}, and \cite[Theorem]{B71} that if 
$\ell^\Phi(\mathbb{Z})$ is reflexive, then its Boyd indices 
(whose definitions can be found in \cite{B71} and also in \cite{KS24}) 
are nontrivial.
So, all hypotheses of \cite[Lemma~6.2]{KT-JAT} are fulfilled.
That lemma yields that if $a\in C[-\pi,\pi]\cap BV[-\pi,\pi]$, then
\[
\|\sigma_n(a)-a\|_{M_\Phi}\to 0
\quad\mbox{as}\quad
n\to\infty,
\]
where the partial sums of the Fourier series of $a$
are defined by
\[
[s_n(a)](\theta):=\sum_{k=-n}^n \widehat{a}_ke^{ik\theta},
\quad
\theta\in\mathbb{R},
\quad
n\in\mathbb{Z}_+,
\]
and the Fej\'er (or Fej\'er-Ces\`aro) means of $a$ are defined by
\[
\sigma_n(a):=\frac{1}{n+1}\sum_{{k=0}}^n s_k(a),
\quad n\in\mathbb{Z}_+.
\]
It is clear that $\sigma_n(a)\in\mathcal{TP}$. So, $a\in C_\Phi$.
This completes the proof of the first inclusion. The second
inclusion follows from \eqref{eq:embedding-of-multipliers}.
\end{proof}
\subsection{Compactness of Hankel operators with symbols in $C_\Phi$}
The following lemma extends \cite[Theorem~2.47(a)]{BS06}.
\begin{lemma}\label{le:Hankel-compactness}
Let $\Phi$ be an $N$-function such that the corresponding Orlicz sequence
space $\ell^\Phi(\mathbb{Z})$ is reflexive. If $a\in C_\Phi$, then the Hankel
operators $H(a)$ and $H(\widetilde{a})$ are compact on the subspace
$\ell^\Phi(\mathbb{Z}_+)$.
\end{lemma}
\begin{proof}
If $a\in C_\Phi$, then there exists a sequence $a_n\in\mathcal{TP}$
such that
\begin{equation}\label{eq:Hankel-compactness-1}
\|L(a-a_n)\|_{\mathcal{B}(\ell^\Phi(\mathbb{Z}))}
=
\|a-a_n\|_{M_\Phi}=o(1)
\quad\mbox{as}\quad
n\to\infty.
\end{equation}
Since $a_n\in\mathcal{TP}$, it follows from 
\eqref{eq:Ha-representation}--\eqref{eq:Ha-tilde-representation}
that $H(a_n)$ and $H(\widetilde{a_n})$ are finite-rank operators.
On the other hand, \eqref{eq:Hankel-compactness-1} yields that
\begin{align*}
\|H(a)-H(a_n)\|_{\mathcal{B}(\ell^\Phi(\mathbb{Z}_+))}
&=
\|PL(a-a_n)QJ\|_{\mathcal{B}(\ell^\Phi(\mathbb{Z}))}
\le 
\|L(a-a_n)\|_{\mathcal{B}(\ell^\Phi(\mathbb{Z}))}
=o(1),
\\
\|
H(\widetilde{a})-H(\widetilde{a_n})
\|_{\mathcal{B}(\ell^\Phi(\mathbb{Z}_+))}
&=
\|JQL(a-a_n)P\|_{\mathcal{B}(\ell^\Phi(\mathbb{Z}))}
\le 
\|L(a-a_n)\|_{\mathcal{B}(\ell^\Phi(\mathbb{Z}))}
=o(1)
\end{align*}
as $n\to\infty$, which imply that $H(a)$ and $H(\widetilde{a})$
are compact as uniform limits of sequences of finite-rank
operators.
\end{proof}
The following corollary of Lemmas~\ref{le:Widom-formulas} 
and~\ref{le:Hankel-compactness}
is analogous to \cite[Lemma~6.3]{BSey00}.
\begin{corollary}\label{co:cosets-commute}
Let $\Phi$ be an $N$-function such that the corresponding Orlicz sequence
space $\ell^\Phi(\mathbb{Z})$ is reflexive. If $a\in M_\Phi$ and {$f\in C_\Phi$,}
then
\[
T^\pi(af) = T^\pi(a) T^\pi(f)=T^\pi(f) T^\pi(a).
\]
\end{corollary}
\subsection{Inclusion $M_{\langle\Phi\rangle}\subset M_\Phi$}
Let us show that $M_{\langle\Phi\rangle}$ is contained in $M_\Phi$.
\begin{lemma}\label{le:inclusion}
Let $\Phi$ be an $N$-function with the Matuszewska-Orlicz indices
satisfying \eqref{eq:indices-nontriviality}. If $\theta$ satisfies
\eqref{eq:small-theta} and $\Phi_\theta$ is defined by 
\eqref{eq:Phi-theta}--\eqref{eq:Psi-theta}, then 
$M_{\Phi_\theta}\subset M_\Phi$.
\end{lemma}
\begin{proof}
Let $a\in M_{\Phi_\theta}$. It follows from Lemma~\ref{le:reflexivity-indices},
Theorem~\ref{th:interpolation-i}, and \eqref{eq:lower-upper-MO-indices}
that the space $\ell^{\Phi_\theta}(\mathbb{Z})$ is reflexive. Inequality
\eqref{eq:embedding-of-multipliers} and the equality
\begin{equation}\label{eq:ell-2-multipliers}
\|L(a)\|_{\mathcal{B}(\ell^2(\mathbb{Z}))}
=
\|a\|_{L^\infty(-\pi,\pi)}
\end{equation}
(see, e.g., \cite[Theorem~1.1 and equalities (1.5)--(1.6)]{BS99} 
and also
\cite[Section~XXIII.2]{GGK93}, \cite[Section~3.1]{GGK03}) imply that
\[
\|L(a)\|_{\mathcal{B}(\ell^2(\mathbb{Z}))}
\le 
\|L(a)\|_{\mathcal{B}(\ell^{\Phi_\theta}(\mathbb{Z}))}.
\]
Combining this inequality with Theorem~\ref{th:interpolation-ii},
we get
\begin{align*}
\|a\|_{M_\Phi}
&=
\|L(a)\|_{\mathcal{B}(\ell^\Phi(\mathbb{Z}))}
\le 
C_{\Phi,\theta}
\|L(a)\|_{\mathcal{B}(\ell^{\Phi_\theta}(\mathbb{Z}))}^{1-\theta}
\|L(a)\|_{\mathcal{B}(\ell^2(\mathbb{Z}))}^\theta
\\
&\le
C_{\Phi,\theta}
\|L(a)\|_{\mathcal{B}(\ell^{\Phi_\theta}(\mathbb{Z}))}
=
C_{\Phi,\theta}
\|a\|_{M_{\Phi_\theta}}.
\end{align*}
Thus $a\in M_\Phi$ and $M_{\Phi_\theta}\subset M_\Phi$.
\end{proof}
Lemma~\ref{le:inclusion} immediately implies that 
$M_{\langle\Phi\rangle}\subset M_\Phi$.
\subsection{Proof of Theorem~\ref{th:main}}
We will follow the arguments presented in \cite[Section~2.69]{BS06}
and \cite[Section~6]{BSey00}. Let
\[
\mathcal{A}:=
\mathcal{B}(\ell^\Phi(\mathbb{Z}_+))/\mathcal{K}(\ell^\Phi(\mathbb{Z}_+))
\]
and $\mathcal{T}:=[-\pi,\pi)$. Let $C^1[-\pi,\pi]$ denote the set of
all continuously differentiable $2\pi$-periodic functions 
$f:\mathbb{R}\to\mathbb{C}$. It is well known that
\begin{equation}\label{eq:proof-main-1}
C^1[-\pi,\pi]\subset C[-\pi,\pi]\cap BV[-\pi,\pi].
\end{equation}
For $\tau\in\mathcal{T}$, define $\mathcal{N}_\tau$ as the set of all
functions $f\in C^1[-\pi,\pi]$ with the following properties:
$0\le f\le 1$;
$f=0$ in $[\tau-\pi,\tau+\pi)\setminus W_{f,\tau}$,
where $W_{f,\tau}$ is an open interval of length less than 
$2\pi$ centered at $\tau$;
$f=1$ in an open interval $U_{f,\tau}$ centered at $\tau$ and such
that $U_{f,\tau}\ne\emptyset$ and $W_{f,\tau}\setminus U_{f,\tau}\ne\emptyset$;
$f$ is increasing in
\[
(W_{f,\tau}\setminus U_{f,\tau})_-:=
\{x\in W_{f,\tau}\setminus U_{f,\tau}\ :\ x<\tau\};
\]
$f$ is decreasing in
\[
(W_{f,\tau}\setminus U_{f,\tau})_+:=
\{x\in W_{f,\tau}\setminus U_{f,\tau}\ :\ \tau<x\}.
\]
It is clear that $\|f\|_{L^\infty(-\pi,\pi)}=1$ and $V(f)=2$ for every
$f\in\mathcal{N}_\tau$. The Stechkin inequality 
(see \cite[Theorem~1.2(c)]{KT-JAT}) implies that $f\in M_\Phi$ and hence
$T^\pi(f)\in\mathcal{A}$ for all $f\in\mathcal{N}_\tau$. Put
\[
\mathcal{M}_\tau:=\{T^\pi(f)\ : \ f\in\mathcal{N}_\tau\},
\quad
\tau\in\mathcal{T}.
\]

Let us show that the sets $\mathcal{M}_\tau$, $\tau\in\mathcal{T}$, have the 
properties required in the Gohberg-Krupnik local principle 
(see Theorem~\ref{th:Gohberg-Krupnik}). Fix $\tau\in\mathcal{T}$.
Since $f\ne 0$ for all $f\in\mathcal{N}_\tau$, it follows from
Corollary~\ref{co:noncompactness-Wiener-Hopf} that $T^\pi(f)\ne 0$,
whence $0\notin\mathcal{M}_\tau$.

Let $T^\pi(f_1), T^\pi(f_2)\in\mathcal{M}_\tau$. Then there exists
$f\in\mathcal{N}_\tau$ with so small support that $f_jf=ff_j=f$ for
$j=1,2$. It follows from this equality, inclusions
\eqref{eq:proof-main-1}, Lemma~\ref{le:CBV-in-C-Phi}
and Corollary~\ref{co:cosets-commute} that
\[
T^\pi(f_j)T^\pi(f)=T^\pi(f)T^\pi(f_j)=T^\pi(f),
\quad
j=1,2.
\]
Now let $\{T^\pi(f_\tau)\}_{\tau\in\mathcal{T}}$ be a collection of elements
$T^\pi(f_\tau)\in\mathcal{M}_\tau$. Since all elements of $\mathcal{N}_\tau$
are $2\pi$-periodic and $[-\pi,\pi]$ is compact, there is a finite subollection 
$\{f_{\tau_j}\}_{j=1}^m$ such that
\[
g:=f_{\tau_1}+\dots+f_{\tau_m}\ge 1.
\]
Hence $g^{-1}$, defined by $g^{-1}(\theta):=1/g(\theta)$ for 
$\theta\in\mathbb{R}$, belongs to $C^1[-\pi,\pi]$. It follows from the equality $g^{-1}g=gg^{-1}=1$, 
inclusion \eqref{eq:proof-main-1}, Lemma~\ref{le:CBV-in-C-Phi},
and Corollary~\ref{co:cosets-commute} that
\[
T^\pi(g^{-1})T^\pi(g)=T^\pi(g)T^\pi(g^{-1})=T^\pi(1)=I^\pi,
\]
which shows that $T^\pi(g)$ is invertible in $\mathcal{A}$.
Similarly, if $a\in M_{\langle\Phi\rangle}\subset M_\Phi$
(see Lemma~\ref{le:inclusion}) and $f\in\mathcal{N}_\tau$
for $\tau\in\mathcal{T}$, then 
\[
T^\pi(a)T^\pi(f)=T^\pi(f)T^\pi(a)
\]
in view of Corollary~\ref{co:cosets-commute}. So 
$T^\pi(a)\in\operatorname{Com}\mathcal{M}$, where
\[
\mathcal{M}:=\bigcup_{\tau\in\mathcal{T}}\mathcal{M}_\tau.
\]

Suppose that $a\in M_{\langle\Phi\rangle}$ and for each $\tau\in\mathcal{T}$
there exists $a_\tau\in M_{[\Phi]}$ such that
\begin{equation}\label{eq:proof-main-2}
\operatorname{dist}_\tau(a,a_\tau)=0
\end{equation}
and $T(a_\tau)$ is Fredholm on $\ell^\Phi(\mathbb{Z}_+)$. Then the coset
$T^\pi(a_\tau)$ is invertible in the Calkin algebra $\mathcal{A}$
(see, e.g., \cite[Ch.~XI, Theorem~5.2]{GGK90}).

Fix $\tau\in\mathcal{T}$. Let us show that $T^\pi(a)$ and $T^\pi(a_\tau)$
are $\mathcal{M}_\tau$-equivalent. Since $a\in M_{\langle\Phi\rangle}$
and $a_\tau\in M_{[\Phi]}$,
there exists $\theta>0$ satisfying \eqref{eq:small-theta} and such that
$a,a_\tau\in M_{\Phi_\theta}$, where $\Phi_\theta$ is given by
\eqref{eq:Phi-theta}--\eqref{eq:Psi-theta}. For every $f\in\mathcal{N}_\tau$,
it follows from inclusion \eqref{eq:proof-main-1}, Lemma~\ref{le:CBV-in-C-Phi},
and Corollary~\ref{co:cosets-commute} that
\[
\left(T^\pi(a)-T^\pi(a_\tau)\right)T^\pi(f)=T^\pi((a-a_\tau)f).
\]
Therefore,
\begin{align}
\inf_{f\in\mathcal{N}_\tau}
\left\|\left(T^\pi(a)-T^\pi(a_\tau)\right)T^\pi(f)\right\|_\mathcal{A}
&=
\inf_{f\in\mathcal{N}_\tau}
\left\|T^\pi((a-a_\tau)f)\right\|_\mathcal{A}
\nonumber\\
&\le 
\inf_{f\in\mathcal{N}_\tau}
\left\|L((a-a_\tau)f)\right\|_{\mathcal{B}(\ell^\Phi(\mathbb{Z}))}.
\label{eq:proof-main-3}
\end{align}
Theorem~\ref{th:interpolation-ii} implies that
\begin{align}
&
\inf_{f\in\mathcal{N}_\tau}
\|L((a-a_\tau)f)\|_{\mathcal{B}(\ell^\Phi(\mathbb{Z}))}
\nonumber\\
&\quad \le 
C_{\Phi,\theta}
\inf_{f\in\mathcal{N}_\tau}
\left(
\|L((a-a_\tau)f)\|_{\mathcal{B}(\ell^{\Phi_\theta}(\mathbb{Z}))}^{1-\theta}
\|L((a-a_\tau)f)\|_{\mathcal{B}(\ell^2(\mathbb{Z}))}^\theta
\right).
\label{eq:proof-main-4}
\end{align}
It follows from Theorem~\ref{th:interpolation-i} and Stechkin's inequality
(see \cite[Theorem~1.2(c)]{KT-JAT}) that for all $f\in\mathcal{N}_\tau$,
\begin{align}
\left\|L((a-a_\tau)f)\right\|_{\mathcal{B}(\ell^{\Phi_\theta}(\mathbb{Z}))}
&\le 
\left\|L(a-a_\tau)\right\|_{\mathcal{B}(\ell^{\Phi_\theta}(\mathbb{Z}))}
\left\|L(f)\right\|_{\mathcal{B}(\ell^{\Phi_\theta}(\mathbb{Z}))}
\nonumber\\
&\le
C_{\ell^{\Phi_\theta}}
\|a-a_\tau\|_{M_{\Phi_\theta}} 
\left(\|f\|_{L^\infty(-\pi,\pi)}+V(f)\right)
\nonumber\\
&=
3C_{\ell^{\Phi_\theta}}\|a-a_\tau\|_{M_{\Phi_\theta}},
\label{eq:proof-main-5}
\end{align}
where $C_{\ell^{\Phi_\theta}}>0$ is the constant in the Stechkin 
inequality. On the other hand, 
\begin{equation}\label{eq:proof-main-6}
\left\|L((a-a_\tau)f)\right\|_{\mathcal{B}(\ell^2(\mathbb{Z}))}
=
\|(a-a_\tau)f\|_{L^\infty(-\pi,\pi)}
\end{equation}
(see \eqref{eq:ell-2-multipliers}).
Combining \eqref{eq:proof-main-3}--\eqref{eq:proof-main-6}, we see that
\begin{align}
&
\inf_{f\in\mathcal{N}_\tau}
\left\|\left(T^\pi(a)-T^\pi(a_\tau)\right)T^\pi(f)\right\|_\mathcal{A}
\nonumber\\
&\quad
\le 
C_{\Phi,\theta}
\left(
3C_{\ell^{\Phi_\theta}} \|a-a_\tau\|_{M_{\Phi_\theta}}
\right)^{1-\theta}
\left(
\inf_{f\in\mathcal{N}_\tau}
\|(a-a_\tau)f\|_{L^\infty(-\pi,\pi)}\right)^\theta.
\label{eq:proof-main-7}
\end{align}

Let $\mathcal{U}_\tau$ be the family of all open finite intervals containing 
the point $\tau\in\mathcal{T}$. If $f\in\mathcal{N}_\tau$, then there
exists $u\in\mathcal{U}_\tau$ such that $f|_u=1$. Hence
\[
\|a|_u-a_\tau|_u\|_{L^\infty(u)}
\le 
\|(a-a_\tau)f\|_{L^\infty(-\pi,\pi), }
\]
which implies that
\begin{equation}\label{eq:proof-main-8}
\operatorname{dist}_\tau(a,a_\tau)
=
\inf_{u\in\mathcal{U}_\tau}
\|a|_u-a_\tau|_u\|_{L^\infty(u)}
\le 
\inf_{f\in\mathcal{N}_\tau}
\|(a-a_\tau)f\|_{L^\infty(-\pi,\pi)}.
\end{equation}
On the other hand, if $u\in\mathcal{U}_\tau$, then there exists 
$f\in\mathcal{N}_\tau$ such that $\operatorname{supp}f\subset u$. Then
\[
\|(a-a_\tau)f\|_{L^\infty(-\pi,\pi)}
\le 
\|a|_u-a_\tau|_u\|_{L^\infty(u)},
\]
which yields that
\begin{equation}\label{eq:proof-main-9} 
\inf_{f\in\mathcal{N}_\tau}
\|(a-a_\tau)f\|_{L^\infty(-\pi,\pi)}
\le 
\inf_{u\in\mathcal{U}_\tau}
\|a|_u-a_\tau|_u\|_{L^\infty(u)}
=
\operatorname{dist}_\tau(a,a_\tau).
\end{equation}
Equality \eqref{eq:proof-main-2} and inequalities
\eqref{eq:proof-main-8}--\eqref{eq:proof-main-9}
imply that
\[
\inf_{f\in\mathcal{N}_\tau}
\|(a-a_\tau)f\|_{L^\infty(-\pi,\pi)}
=
\operatorname{dist}_\tau(a,a_\tau)=0.
\]
It follows from this equality and \eqref{eq:proof-main-7}
that
\[
\inf_{f\in\mathcal{N}_\tau}
\left\|\left(T^\pi(a)-T^\pi(a_\tau)\right)T^\pi(f)\right\|_\mathcal{A}=0.
\]
Analogously, one can prove that
\[
\inf_{f\in\mathcal{N}_\tau}
\left\|T^\pi(f)\left(T^\pi(a)-T^\pi(a_\tau)\right)\right\|_\mathcal{A}=0.
\]
Thus $T^\pi(a)$ and $T^\pi(a_\tau)$ are $\mathcal{M}_\tau$-equivalent.
Applying the Gohberg-Krupnik local principle 
(see Theorem~\ref{th:Gohberg-Krupnik}), we conclude that the coset
$T^\pi(a)$ is invertible in the Calkin algebra $\mathcal{A}$. 
Therefore, the discrete Wiener-Hopf operator $T(a)$ is Fredholm
on the space $\ell^\Phi(\mathbb{Z})$ in view of 
\cite[Ch.~XI, Theorem~5.2]{GGK90}.
\qed
\subsection{A slightly more general version of Theorem~\ref{th:main}}
A careful analysis of the proof of Theorem~\ref{th:main} presented in the 
previous subsection shows that the following slightly more general version 
of it is true.
\begin{theorem}
Let $\Phi$ be an $N$-function with the Matuszewska-Orlicz indices satisfying
\eqref{eq:indices-nontriviality}. Suppose that $a\in M_{\langle\Phi\rangle}$
and for every $\tau\in[-\pi,\pi)$ there exist $\theta(\tau)$ and $a_\tau$ 
such that
\begin{enumerate}
\item[{\rm(i)}]
$0<\theta(\tau)<2\min\{1/\beta_\Phi,1-1/\alpha_\Phi\}$;

\item[{\rm(ii)}]
$a,a_\tau\in M_{\Phi_{\theta(\tau)}}$, where $\Phi_{\theta(\tau)}$
is defined by \eqref{eq:Phi-theta}--\eqref{eq:Psi-theta} with $\theta$
replaced by $\theta(\tau)$;

\item[{\rm(iii)}]
$\operatorname{dist}_\tau(a,a_\tau)=0$;

\item[{\rm(iv)}]
the discrete Wiener-Hopf operator $T(a_\tau)$ is Fredholm on
$\ell^\Phi(\mathbb{Z}_+)$.
\end{enumerate}
Then the discrete Wiener-Hopf operator $T(a)$ is Fredholm on 
$\ell^\Phi(\mathbb{Z}_+)$.
\end{theorem}
\section*{Declarations}
\subsection*{Acknowledgements}
The first author would like to thank Professor Mieczys{\l}aw Masty{\l}o
for useful discussions concerning 
Theorems~\ref{th:interpolation-i} and~\ref{th:interpolation-ii}
during the ``Workshop on Banach spaces and Banach lattices", 
Madrid, Spain, September 9--13, 2019.
We would like to thank Professor Albrecht B\"ottcher and anonymous
referees for useful remarks.

\subsection*{Funding}
This work is funded by national funds through the FCT - Funda\c{c}\~ao para a 
Ci\^encia e a Tecnologia, I.P., under the scope of the projects 
UIDB/00297/2020 (\url{https://doi.org/10.54499/UIDB/00297/2020})
and 
UIDP/ 00297/2020
(\url{https://doi.org/10.54499/UIDP/00297/2020})
(Center for Mathematics and Applications).
The second author is funded by national funds through the FCT – 
Funda\c{c}\~ao para a Ci\^encia e a Tecnologia, I.P., 
under the scope of the PhD scholarship UI/BD/152570/2022.
\subsection*{Conflict of interest}
The authors declare no competing interests.
\subsection*{Data availability} 
The manuscript does not contain any associated data.
\bibliography{OKST-Grudsky70}
\end{document}